\documentclass[12pt]{article}
\usepackage{amssymb,amsmath}
\newtheorem{theorem}{Theorem}[section]
\newtheorem{lemma}{Lemma}[section]
\newtheorem{proposition}{Proposition}[section]
\newtheorem{corollary}{Corollary}[section]
\newtheorem{definition}{Definition}[section]

\begin{document}

$\;$

\begin{center}

\LARGE{Compact special Legendrian surfaces in $S^5$}

\vspace{1pc}

\large{Sung Ho Wang} \\

\vspace{0.7pc}

 Department of Mathematics \\
 Postech \\
 Pohang, Korea 790 - 784 \\
 \emph{email}: \textbf{wang@postech.ac.kr}
\end{center}

\vspace{2pc}

\thispagestyle{empty}

\begin{center}
\textbf{Abstract}
\end{center}

A surface $\Sigma \subset S^5 \subset \mathbb{C}^3$ is called
\emph{special Legendrian} if the cone $0 \times \Sigma \subset
\mathbb{C}^3$ is special Lagrangian. The purpose of this paper is
to propose a general method toward constructing compact special Legendrian
surfaces of high genus. It is proved \emph{there exists a compact,
orientable, Hamiltonian stationary Lagrangian surface of genus
$1+\frac{k(k-3)}{2}$ in $\, \mathbb{C}P^2$ for each integer $\, k \geq 3$,
which is a smooth branched surface except at most finitely many conical singularities.}
If this surface is smooth, it is minimal and the Legendrian lift of the surface
is the desired compact special Legendrian surface.

We first establish the existence of a minimizer of area
among Lagrangian disks in a relative homotopy class of a
K\"ahler-Einstein surface without Lagrangian homotopy classes with
respect to a configuration $\, \Gamma$ that consists of the fixed
point loci of K\"ahler involutions. $\, \Gamma$ in addition must satisfy
certain null relative homotopy conditions and angle
criteria. The fundamental domain thus obtained
is smooth along the boundary, and has finitely many interior singular points.
We then apply successive \emph{reflection} of this
fundamental domain along its boundary to obtain a complete or
compact Lagrangian surface.

\vspace{1pc}


\noindent \textbf{MS classification}: 53C25

\thispagestyle{empty}

\newpage
\setcounter{page}{1}

$\;$

\begin{center}
\LARGE{ Compact special Legendrian surfaces in $S^5$ }

\vspace{1pc}

\large{Sung Ho Wang}

\end{center}

\vspace{2pc}

\begin{center}  \textbf{Table of Contents} \end{center}

$\qquad \quad$ 0. Introduction

$\qquad \quad$ 1. Special Lagrangian cones in $\mathbb{C}^{n+1}$

$\qquad \quad$ 2. Special Legendrian surfaces in $S^5$

$\qquad \quad$ 3. Reflection principle

$\qquad \quad$ 4. Rational geodesic polygon

$\qquad \quad$ 5. Construction procedure

$\qquad \quad$ 6. Surfaces $\, \Sigma_{k, 3}$

$\qquad \quad$ 7. Gauss map, polar surface, and bipolar surface

$\qquad \quad$ 8. Intrinsic characterization and dual reflection
principle

\vspace{2pc}

\begin{center}
\noindent \textbf{\large{Introduction}}
\end{center}

Bryant showed every compact Riemann surface admits conformal minimal
embeddings in $S^4$ with vanishing complex quartic form [Br1]. With
regard to the fibration
\begin{equation}
Sp(2) \to \mathbb{C}P^3 \to S^4 = \mathbb{H}P^1, \notag
\end{equation}
such an embedding corresponds to a holomorphic Legendrian curve in
$\mathbb{C}P^3$ under $Sp(2)$ invariant holomorphic contact structure.
Similarly, every compact Riemann surface can be immersed as a branched
complex curve in $S^6$ with null torsion with respect to $G_2$ invariant
almost complex structure via fibration
\begin{equation}
 G_2/U(2) = Gr^+(2, \mathbb{R}^7) \to S^6. \notag
\end{equation}
The lift of such a complex curve is integral to  a rank 3 holomorphic
Pfaffian system on $Gr^+(2, \mathbb{R}^7) \subset \mathbb{C}P^6$, for
which Cartan gave a local normal form in terms of
birational coordinates [Br2]. These two differential systems for minimal
immersion of surfaces have the same local generality of solutions,
i. e., one holomorphic function of one complex variable [BCG].

The twistor constructions above, as they stand, do not apply to reduce
to holomorphic category the differential systems for minimal immersion
of surfaces in odd dimensional spheres with the same local generality
of solutions. For instance, minimal surfaces in $S^3$ correspond to
complex curves in $SO(4)/SO(2)$ with respect to a nonintegrable
$CR$-structure, Section 7.

A surface $\Sigma \subset S^5 \subset \mathbb{C}^3$ is called
\emph{special Legendrian} if the cone $0 \times \Sigma \subset \mathbb{C}^3$ is
special Lagrangian. Up to $U(3)$ motion, a special Legendrian surface is
a minimal Legendrian surface, and hence corresponds to a minimal Lagrangian
surface in $\mathbb{C}P^2$ via Hopf map.

The purpose of this paper is to propose a general method toward constructing
compact special Legendrian surfaces of high genus. It is proved
\emph{there exists a compact, orientable, Hamiltonian stationary
Lagrangian surface of genus $1+\frac{k(k-3)}{2}$ in $\mathbb{C}P^2$ for
each integer $k \geq 3$, which is a smooth branched surface except at most finitely
many conical singularities.}
If this surface is smooth, it is minimal and the Legendrian lift of the surface,
Proposition 1.6, is then the desired compact special Legendrian surface.
As of this writing, we are not able to resolve the singularities.

Our method of construction is analogous to that of Lawson's on compact
minimal surfaces in $S^3$ [La1]. We first find a minimizer of certain
free boundary problem, and apply successive \emph{reflection} of this fundamental
domain to obtain a complete Lagrangian surface. If the associated
reflection group is finite, the resulting surface is compact.
The reflection principle, to be explained below, is motivated from
the reflection principle for minimal surfaces in $S^3$, and the similarity
between the elliptic Monge-Ampere systems describing special Legendrian surfaces
and minimal surfaces in $S^3$.

Let $\phi$ be a calibration on a Riemannian manifold $X$. Consider a compact
calibrated submanifold $M$ with boundary $\partial M \subset N$, where $N$ is
a $\phi$-null submanifold, a submanifold on which $\phi=0$. Then via Stokes
theorem, $M$ is absolutely mass minimizing in its relative homology class
with respect to $N$. Thus $\phi$-null submanifolds serve as natural
constraints for free boundary problem for $\phi$-submanifolds.
For example, finding a special Lagrangian submanifold with boundary
in a given complex hypersurface in a Calabi-Yau manifold roughly corresponds
to solving Neumann problem.

Since such $\phi$-calibrated submanifold intersects a $\phi$-null submanifold
orthogonally along the boundary, the observation above leads to the following
reflection principles;

\emph{1. Suppose $N$ is a $\phi$-null submanifold
which is a fixed point locus of a $\, \pm \phi$ isometric involution
$\, \sigma$ of $\, X$, i. e., $\, \sigma^* \phi = \pm \phi$.
Then any $\phi$-submanifold with boundary in $N$ and $C^{1}$ up to
boundary can be analytically continued across $N$ via $\sigma$.}

\emph{ 2. Suppose $M$ is a $\phi$-calibrated submanifold which is a fixed point
locus of $\, \phi$ isometric involution $\, \sigma$. Then any minimal $\phi$-null
submanifold with boundary in $M$ and $C^1$ up to boundary can be analytically
continued across $M$ via $\sigma$.}

Basic facts and fundamental equations for special Lagrangian
submanifolds and special Legendrian surfaces are recorded in Section 1
and 2. After introducing the reflection principle for special
Lagrangian 3-folds with complex, or anti-special Lagrangian boundary,
we describe a set of conditions for the constraining polygon that is
necessary for our construction. Section 5 contains the main theorem
that there exists a minimizer of  area of the free boundary problem for
Lagrangian surfaces in $\mathbb{C}P^{2}$, which is smooth along the
boundary and with finitely many interior singular points,
if the constraining polygon satisfy certain vanishing relative homotopy
conditions. It is based upon the fundamental work of Schoen and Wolfson [ScW],
and holds for general K\"ahler-Einstein surfaces without Lagrangian
homotopy classes. The main theorem is then used in Section 6 to
construct a family of compact Hamiltonian stationary Legendrian surfaces. The
constraining variety in this case is a complex geodesic triangle
derived from the generating set of a finite unitary reflection
group in $\mathbb{C}^{3}$ [Cox]. Section 7 is devoted to the study
of surfaces associated to a given special Legendrian surface. We show in
particular that a compact minimal Lagrangian surface in $\mathbb{C}P^2$
cannot lie in an open geodesic ball of radius
$\frac{\pi}{2} - \arccos(\frac{1}{\sqrt{3}})$.
In Section 8, we give an intrinsic characterization of the induced
Riemannian metrics on special Legendrian surfaces or minimal
Lagrangian surfaces in complex space forms.

We'd like to thank Nancy for her continuing support.

\section{ Special Lagrangian cones in $\mathbb{C}^{n + 1}$}

Let $\mathbb{C}^{n+1}$ be the complex vector space of dimension
$n+1$, which we always regard as the real vector space $\mathbb{R}^{2n+2}$
with coordinates $(x, y)$ and the complex structure
\begin{displaymath}
J( x, y ) = ( -y, x ),
\end{displaymath}
where  $x = ( \, x^1, \, . . . \,  x^{n+1} \, ),
y = ( \, y^1, \, . . . \,  y^{n +1} \, )$.
Let
\begin{align}
\varpi &= \sum_k dx^k \wedge dy^k \notag \\
\langle v, w \rangle &= \varpi(v, Jw) \quad \,  \mbox{for} \, \,  \; v, w \, \in
\mathbb{C}^{n+1} \notag
\end{align}
be the standard K\"ahler form and the metric. Set $z^k = x^k + i \, y^k$,
and we denote
\begin{align}
\Upsilon &= dz^1 \wedge \; . . . \;  dz^{n+1} \notag \\
\Upsilon_{\theta} &= e^{i \theta} \, \Upsilon \notag
\end{align}
the holomorphic volume form, and holomorphic volume form of phase $e^{i \theta}$.

Let $Gr^{+}_{\mathbb{R}}(k, \mathbb{C}^{n+1})$ be the set of
oriented real Grassmann $k$-planes in $\mathbb{C}^{n+1}$. An
element $\sigma \in Gr^{+}_{\mathbb{R}}(k, \mathbb{C}^{n+1})$ is
called \emph{isotropic} if $\varpi|_{\sigma}=0$, $1 \leq k \leq
n+1 \,$ necessarily, and it is called \emph{Lagrangian} when
$k=n+1$. We denote
\begin{align}
Isot^{+}(k, \mathbb{C}^{n+1}) &= \{ \, \sigma \in
Gr^+_{\mathbb{R}}(k,
\mathbb{C}^{n+1}) \, \, | \, \, \;
                              \varpi|_{\sigma} = 0 \, \} \notag \\
                         &= SU(n+1)/(SO(k) \times SU(n+1-k)) \notag \\
Slag(n+1) &= \{ \; L \in Gr^{+}_{\mathbb{R}}(n+1,
\mathbb{C}^{n+1})
             \; | \; dvol_L = Re ( \Upsilon ) |_{L} \; \} \notag \\
          &= SU(n+1)/SO(n+1)  \notag \\
Slag(n+1)^{\perp}
          &= \{ \; L' \in Gr^{+}_{\mathbb{R}}(n+1,
\mathbb{C}^{n+1})
             \; | \; dvol_{L'} = - Im ( \Upsilon ) |_{L'} \; \}
\notag \\
          &= \{ \, (-1)^{\frac{n}{2}+1} J(L) = L^{\perp} \, | \; L
\in Slag(n+1) \; \} \; \, \mbox{ if $ n+1$ is odd.} \notag
\end{align}
An element of $Slag(n+1)$ is called a \emph{special Lagrangian}
$(n+1)$-plane, and it is characterized, up to choice of right
orientation, as a Lagrangian $(n+1)$-plane on which
$Im( \Upsilon) = 0$.
\begin{lemma}\textnormal{[HL1]}
Given $\sigma \in Isot^{+}(n, \mathbb{C}^{n+1})$, there exists a unique
special Lagrangian $(n+1)$-plane $L_{\sigma}$ such that $\sigma
\subset L_{\sigma}$.
\end{lemma}

Lemma 1.1 provides a double fibration

\begin{picture}(300,60)(-37,0)
\put(130,40){$Isot^{+}(n, \mathbb{C}^{n+1})$}
\put(170,25){$\searrow$}
\put(135,25){$\swarrow$}
\put(119,10){$S^{2n+1}$}
\put(169,10){$Slag(n+1).$}
\put(182,27){$\pi_2$}
\put(123,27){$\pi_1$}
\end{picture}

\noindent Each fiber of the projection $\pi_1$ is identified with $Slag(n)$ and
each fiber of $\pi_2$ with $S^n$.

Let $M$ be an oriented manifold of dimension $n+1$.
\begin{definition}
An immersion $u : M \to \mathbb{C}^{n+1}$ is \textnormal{special Lagrangian} if
\begin{equation}
du(T_p M) \in Slag(n+1) \notag
\end{equation}
for each $p \in M$. The image $u(M)$ of such an immersion is a \textnormal{special
Lagrangian submanifold}.
\end{definition}
Note that a special Lagrangian submanifold is orientable.

\begin{proposition}\textnormal{[HL1]}
An immersion $u :  M \to \mathbb{C}^{n+1}$ is minimal Lagrangian if and only if
\begin{equation}
  u^* \, Re ( \Upsilon_{\theta} ) =  dvol_M \notag
\end{equation}
for some constant $\theta$.
\end{proposition}
Such $ \theta$ is called the \emph{phase} of the minimal Lagrangian submanifold.
Minimal Lagrangian submanifolds of phase $\theta$ are thus integral manifolds
of the differential system
\begin{equation}
 \mathcal{I}_{\theta} = \{ \; \varpi, \, Im ( \Upsilon_{\theta} ) \; \}.
\end{equation}

Since each $Re ( \Upsilon_{\theta} )$ gives a calibration in $\mathbb{C}^{n+1}$,
a minimal Lagrangian submanifold is absolutely mass minimizing in the
relative homology class with respect to its boundary. In particular,
$C^1$ special Lagrangian submanifolds  are real analytic.

The following is a direct consequence of Lemma 1.1 and the fact the
differential system (1) is involutive [BCG].
\begin{proposition}
Let $N^n \subset \mathbb{C}^{n+1}$ be a connected,  analytic, isotropic submanifold
of dimension $n$. There exists a unique maximal connected special Lagrangian
submanifold $M^{n+1}$ that contains $N^n$.
\end{proposition}
We mention that real analyticity assumption on $N$ is necessary
both for existence and uniqueness.

Let $S^{2n+1} = \, \{ \, v \in \mathbb{C}^{n+1} \, | \; \; |v| = 1 \, \}$
with the $U(n+1)$ invariant contact structure generated by the contact 1-form
\begin{equation}
 \vec r \, \lrcorner \, \varpi |_{S^{2n+1}},
\end{equation}
where $\vec r = \sum_k x^k \frac{\partial}{\partial x^k}
+ y^k\frac{\partial}{\partial y^k}$.
A $\,k$-dimensional submanifold of $S^{2n+1}$ is called \emph{isotropic} if it is
tangent to the contact hyperplane distribution, $1 \leq k \leq n$ necessarily,
and it is called \emph{Legendrian} when $\,k = n$.

Let $\Sigma \, $ be an oriented manifold of dimension $n$.
\begin{definition}
An immersion $u : \Sigma \to S^{2n+1} \subset \mathbb{C}^{n+1}$
is  \textnormal{special Legendrian} if the cone $ \, 0 \times u (\Sigma) $
is special Lagrangian. The image $u(\Sigma)$ of such an immersion is a
\textnormal{special Legendrian submanifold}.
\end{definition}
Here the cone is given the product orientation of
$\Sigma \times \mathbb{R^{+}}$.
A special Legendrian submanifold is necessarily orientable.

\begin{proposition}
Let $u : \Sigma \to S^{2n+1}$ be a  minimal Legendrian immersion.

1. The cone $ 0 \times u(\Sigma) $ is minimal Lagrangian in $\mathbb{C}^{n+1}$.
In particular, $C^1$ special Legendrian submanifolds are real analytic.

2. When coupled with the Hopf map $\, \pi : S^{2n+1} \to \mathbb{C}P^n$,
$\,  \pi \circ u : \Sigma \to \mathbb{C}P^n$ is a minimal Lagrangian immersion
with respect to the standard K\"ahler structure of $\, \mathbb{C}P^n$.
\end{proposition}

By definition, special Legendrian submanifolds are integral manifolds
of the involutive differential system
\begin{equation}
\mathcal{I} = \{ \, \vec r \, \lrcorner \, \varpi, \, \varpi,
                    \, \vec r \, \lrcorner \, Im(\Upsilon) \, \}.
\end{equation}
\begin{proposition}
Let $\gamma^{n-1} \subset S^{2n+1}$ be a connected, analytic,
isotropic submanifold of dimension \mbox{$n-1$}.
There exists a unique maximal connected  special
Legendrian submanifold $\, \Sigma^n$ that contains $\gamma^{n-1}$.
\end{proposition}

\paragraph{Example 1.1}
\

1. Let $L \subset \mathbb{C}^{n+1}$ be a special Lagrangian $(n+1)$-plane.
The geodesic sphere $ L \cap S^{2n+1} $ is special Legendrian in $S^{2n+1}$.
Its image under Hopf map is the standard
\begin{equation}
RP^n \subset \mathbb{C}P^n.
\end{equation}

2. Let $M(n, \, \mathbb{C}) = \mathbb{C}^{n^2}$ be the space of $n$-by-$n$ complex
matrices. Then the cone $0 \times SU(n) \subset M( n, \, \mathbb{C})$ is  minimal
Lagrangian  with the associated minimal Lagrangian embedding
\begin{equation}
SU(n)/Z_n \to \mathbb{C}P^{n^2-1},
\end{equation}
where $Z_n$ is the center of $SU(n)$.

3. Let $Sym( n, \, \mathbb{C} ) = \mathbb{C}^{ \frac{n(n+1)}{2}}$ be
the space of symmetric complex $n$-by-$n$ matrices.
Consider $0 \times Slag(n) \subset Sym(n, \mathbb{C})$, where
\begin{align}
Slag(n) &= \{ \, U U^t \in Sym(n, \mathbb{C} ) \, | \; U \in SU(n) \;
\}. \notag
\end{align}
It is a minimal Lagrangian cone with the associated minimal
Lagrangian immersion
\begin{equation}
Slag(n) \to \mathbb{C}P^{\frac{n(n+1)}{2}-1}.
\end{equation}
Note that this observation implies a natural map
\begin{equation}
Slag(n) \to Slag( \frac{n(n+1)}{2} ), \notag
\end{equation}
which is also the restriction of the Veronese embedding
$\mathbb{C}P^{m} \to \mathbb{C}P^{\frac{m(m+1)}{2} -1}$,
$m=\frac{n(n+1)}{2} -1$.

4. Let $A( 2n, \, \mathbb{C} ) = \mathbb{C}^{n(2n-1)}$ be the space of skew
symmetric complex $2n$-by-$2n$ matrices. $U \in SU(2n)$ acts on this space by
$\; A \to U \, A \, U^t$ for $A \in A( 2n, \, \mathbb{C})$. The orbit of the matrix
$\begin{pmatrix}
0 & \textrm{I}_n \\
-\textrm{I}_n & 0
\end{pmatrix}$
is $\, SU(2n) / Sp(n)$,
and the cone $0 \times SU(2n)/Sp(n) \subset A(2n, \, \mathbb{C} )$
is minimal Lagrangian. The associated  minimal Lagrangian immersion is
\begin{equation}
SU(2n)/Sp(n) \to \mathbb{C}P^{(2n+1)(n-1)}.
\end{equation}

5. Let $V^{27}$ be the real vector space of Octonianic Hermitian
$3$-by-$3$ matrices and $V_{\mathbb{C}}$ its complexification. The compact
simple Lie group $E_6$ has a faithful special unitary representation on
$V_{\mathbb{C}}$ with $F_4$ as the stabilizer of the identity matrix. The cone
$0 \times E_6 / F_4$ is minimal Lagrangian in $V_{\mathbb{C}}$ with the
associated minimal Lagrangian immersion
\begin{equation}
E_6/F_4 \to \mathbb{C}P^{26}.
\end{equation}

Examples (4), (5), (6), (7), and (8) are in fact the only
irreducible minimal Lagrangian submanifolds
of complex projective spaces with parallel second fundamental form
[Na].

\paragraph{Example 1.2} Consider the minimal Lagrangian cone
\begin{equation}
\{ \; (z^1, \, . . . \, z^{n+1}) \in \mathbb{C}^{n+1} \; |
 \, \, \;  |z_1| = \, . . . \, |z_{n+1}|, \; Im(z_1 z_2 \, . . . \,
z_{n+1}) = 0 \; \}. \notag
\end{equation}
It has the hexagonal $n$-torus as its link, where a hexagonal
$n$-torus is the quotient of Clifford $(n+1)$-torus by the diagonal
action of $S^1$. Hopf map induces an $(n+1)$-fold cover onto its image,
which is again called a hexagonal torus in $\mathbb{C}P^n$.

Special Legendrian differential system on $S^{2n+1}$ can be considered as a
compact real form of the differential system describing complex affine spheres in
$\mathbb{C}^{n+1}$. In this perspective, the product of minimal Lagrangian
submanifolds below is a compact real form of the product for affine spheres
introduced by Calabi [Cal].

\begin{proposition}
Let $M^{n_i}_i \subset \mathbb{C}P^{n_i}, n_i \geq 0,$ be a minimal
Lagrangian submanifold for $i = 1, \, . . . k$. There exists
a product minimal Lagrangian submanifold
\begin{equation}
( M_1 \, \otimes \, . . . \, \otimes \, M_k )^n \, \subset
\, \mathbb{C}P^n \notag
\end{equation}
where $n = \sum_i (n_i + 1) - 1$.
\end{proposition}
\noindent \emph{Proof.} $\; $ Let $\pi^{-1}(M_i) \subset S^{2n_i+1}$
be the inverse image of $M_i$ under the Hopf map, which is a minimal submanifold.
We scale each of $\pi^{-1}(M_i)$ by a factor $r_i=(\frac{n_i+1}{n+1})^{\frac{1}{2}}$,
denote it $\tilde{M}_i$, so that the product
$  \tilde{M}_1 \, \times \, . . . \, \times
\,\tilde{M}_k  \, \subset  S^{n+1}$ is a minimal submanifold.
Now $\pi ( \tilde{M}_1 \, \times \, . . . \, \times  \,\tilde{M}_k )
\, \subset \mathbb{C}P^n$ is  a minimal Lagrangian submanifold. $\square$

Note that
\begin{equation}
 M_1 \, \otimes \, . . . \, \otimes \, M_k  \to
 M_1 \, \times \, . . . \, \times \, M_k \notag
\end{equation}
is a $(k-1)$ dimensional torus bundle. In fact,
minimal Lagrangian torus in $\mathbb{C}P^n$ of Example 1.2 is
a product of $n+1$ minimal Lagrangian submanifolds $\mathbb{R}P^0
\subset \mathbb{C}P^0$.
We also mention that the product of Hamiltonian stable compact minimal
Lagrangian submanifolds is again Hamiltonian stable, Theorem 7.3, [Oh].

In view of Proposition 1.5 and [Na], we have a complete description
of minimal Lagrangian submanifolds in complex projective spaces with
parallel second fundamental form.
\begin{theorem}
Let $M^n \subset \mathbb{C}P^n$ be a minimal Lagrangian submanifold
with parallel second fundamental form. Then it is a part of a product of
minimal Lagrangian submanifolds (4), (5), (6), (7), and (8).
\end{theorem}

Theorem 1.1 can be applied to improve the existing
extrinsic or intrinsic curvature pinching results for minimal
Lagrangian submanifolds [Ma][MRU].

\paragraph{Example 1.3}
Haskins provides families of ruled special Legendrian tori in $S^5$,
and Joyce considers various symmetry reductions [Ha][Jo].

\paragraph{Example 1.4}
Every compact Riemann surface can be embedded in $\mathbb{C}P^3$ as
a holomorphic Legendrian curve with respect to the $Sp(2)$ invariant
holomorphic contact structure [Br1].
Consider the inverse image of a holomorphic Legendrian curve in $S^7
\subset \mathbb{C}^4$ via  Hopf map. Its projection in $\mathbb{C}P^3$
with respect to a different orthogonal complex structure on
$\mathbb{C}^4 = \mathbb{H}^2$ is  minimal Lagrangian.

\begin{proposition} \textnormal{[Wa1]}
Let $u : M \to \mathbb{C}P^n$ be a connected minimal Lagrangian
submanifold, and let
\begin{equation}
\pi_1(M) \to  Hol(M, S^1) \subset O(2) \notag
\end{equation}
denote the holonomy of the associated flat $S^1$ bundle
$u^*(S^{2n+1}) \to M$. Then
\begin{align}
Hol(M, S^1) &\subset \mathbb{Z}_{n+1}  \subset SO(2) \textnormal{\;
if and only if $M$ is
orientable} \notag \\
Hol(M, S^1) &\subset \mathbb{D}_{n+1}  \subset O(2) \textnormal{\; if
$M$ is
nonorientable} \notag
\end{align}
where $\mathbb{D}_{n+1}$ is the dihedral group of order $2(n+1)$.
In particular, if $M$ is compact or embedded, there exists a
connected minimal Legendrian lift $\tilde{M} \subset S^{2n+1}$ that
is compact or embedded respectively.
\end{proposition}

Examples 1.1, 1.2, 1.3, 1.4  and Proposition 1.5 and 1.6
provide many compact special Legendrian (minimal Lagrangian
respectively) submanifolds in $S^{2n+1}$ ($\mathbb{C}P^n$)
with nontrivial topology.

\begin{corollary}
Let $\Sigma^n \subset S^{2n+1}$ be a connected minimal Legendrian submanifold
of even dimension. Then the image of $\, \Sigma$ in $\, \mathbb{C}P^n$ under Hopf
map is nonorientable whenever $\, \Sigma \subset S^{2n+1}$ is invariant under
antipodal involution.
\end{corollary}

\begin{proposition}
Let $\Sigma^n \subset S^{2n+1} $ be a compact special Legendrian
submanifold of positive mass. Let $c(n+1)$ be the isoperimetric constant for
$(n+1)$-dimensional varieties in $\mathbb{R}^{2n+2}$. Then
\begin{equation}
 ||\Sigma|| \geq \frac{c(n+1)}{(n+1)^n}. \notag
\end{equation}
\end{proposition}

\begin{corollary}
There exists a uniform lower bound for the area of compact minimal Lagrangian
submanifolds in $\mathbb{C}P^n$ depending only on the dimension $\, n$.
\end{corollary}
\emph{Proof of the Proposition.}$\; $ Let
$0 \times_1 \Sigma = (0 \times \Sigma) \cap B(1, \mathbb{R}^{2n+2})$, where
$B(1, \mathbb{R}^{2n+2})$ is the unit ball in $\mathbb{R}^{2n+2}$. Since
$0 \times_1 \Sigma$ is mass minimizing,
\begin{equation}
||\Sigma||^{n+1} \geq c(n+1) ||0 \times_1 \Sigma||^n
                  = c(n+1) \frac{||\Sigma||^n}{(n+1)^n}.
           \, \;  \square \notag
\end{equation}

We record the following topological obstruction for Lagrangian embedding.
\begin{proposition}
Let $M^n$ be a compact, orientable, embedded Lagrangian submanifold
in $\mathbb{C}P^n$. Then the Euler Characteristic $\chi(M)=0$.
\end{proposition}
\emph{Proof.} $\; $ Since $M$ is Lagrangian, the tangent bundle $TM$
is isomorphic to the normal bundle $NM$. $M$ is orientable and
embedded, which implies the Euler class of $NM$ comes from the
restriction of an element in $H^n(\, \mathbb{C}P^n, \mathbb{R} ) $
 [Mil]. But $H^*(\, \mathbb{C}P^n, \mathbb{R} )$ is generated by
the K\"ahler form. $\square$

\section{Special Legendrian surfaces in $S^5$}
Let $\Sigma^2 \subset S^5 \subset \mathbb{C}^3$ be a special
Legendrian surface. $\Sigma$  is by definition integral to the
differential system
\begin{equation}
\mathcal{I} = \{ \, \vec r \, \lrcorner \, \varpi, \, \varpi,
                    \, \vec r \, \lrcorner \, Im(\Upsilon) \, \}
\end{equation}
with
\begin{align}
    \varpi &= dx^1\wedge dy^1 + dx^2\wedge dy^2 +
    dx^3\wedge dy^3 \notag \\
- Im(\Upsilon) &= dy^{123} - dy^1\wedge dx^{23} - dy^2 \wedge dx^{31}
                 - dy^3\wedge dx^{12}, \notag
\end{align}
where $dx^{23} = dx^2 \wedge dx^3$, etc.
Let $e_3 : \Sigma \to S^5 \subset \mathbb{R}^6 = \mathbb{C}^3$
be the position function. Take an orthonormal basis $\{  e_1,
e_2 \}$ of $T_{e_3} \Sigma$, which we may regard as an $SO(2)$
equivariant $\mathbb{R}^6$ valued functions on the unit tangent
bundle $S\Sigma$.

Set
\begin{align}
n_i &= J ( e_i ),  \notag \\
e &= ( \, e_1, \, e_2, \, e_3 \, ), \; \;
n= ( \, n_1, \, n_2, \, n_3 \, ). \notag
\end{align}
Since $\Sigma$ is special Legendrian, $e_1 \wedge e_2 \wedge e_3 \, \in
Slag(3)$ and we get an $SO(2)$ equivariant map
\begin{equation}
( \, e, \, n \, ) \, : \, S\Sigma \, \to \, SU(3) \subset SO(6).
\end{equation}
Differentiating (10), we get
\begin{equation}
d (e, n) = (e, n)
\begin{pmatrix}
\alpha & - \beta \\
\beta & \alpha
\end{pmatrix}
\end{equation}
where $\alpha$ is an $\mathfrak{so}(3)$ valued, and $\beta$ is a
trace free $Sym(3, \mathbb{R})$ valued one form.

Let $\{  \omega^1, \, \omega^2  \}$ be the basis of one forms
dual to $\{  e_1, e_2  \}$, and $\rho$ be the connection
form on $S\Sigma$. From the definition,
\begin{equation}
d e_3 = \omega^1 e_1 + \omega^2 e_2  \notag
\end{equation}
and (11)  becomes
\begin{equation}
\alpha = \,
\begin{pmatrix}
0 & \rho &   \omega^1 \\
-\rho & 0 & \omega^2 \\
-\omega^1 & -\omega^2 & 0
\end{pmatrix}, \; \;
\beta = \,
\begin{pmatrix}
\beta^1_1 & \beta^1_2 &  0 \\
\beta^2_1 &  -\beta^1_1  &   0  \\
0 &  0 &  0
\end{pmatrix}.
\end{equation}
Taking the exterior derivative of (11), we obtain the following structure equations.
\begin{align}
d(\omega^1 + i \omega^2) &= i \rho \wedge (\omega^1 + i \omega^2) \\
\beta^1_1 - i \, \beta^1_2 &= \, ( a - i \, b) (\omega^1 + i
\omega^2) \notag \\
d \rho &= K \, \omega^1 \wedge \omega^2 \notag \\
K &= 1 - 2 ( \, a^2 + b^2 \, ). \notag
\end{align}
Put $h = a - i \, b$, then
\begin{align}
dh &= - 3 \, i \, h  \, \rho + h_1 ( \omega^1 + i \omega^2 ) \notag \\
d h_1  &\equiv -4 \, i \, h_1 \, \rho + \frac{3}{2} h \, K \,
(\omega^1 - i \omega^2)
      \; \;\mod \; \; \omega^1 + i \, \omega^2. \notag
\end{align}

Set
\begin{equation}
\Phi = ( \, a - i \, b \, ) \, ( \, \omega^1 + i \,  \omega^2 \, )^3. \notag
\end{equation}
Then $\Phi$ is a well defined cubic differential on $\Sigma$ holomorphic with
respect to the complex structure induced from the metric.

Bonnet type fundamental theorem for special Legendrian surfaces
can now be stated as follows. [Gr].
\begin{theorem}
Consider a  triple $( \, \Sigma, \, \mathfrak{g}, \, \Phi \, )$ of a Riemann
surface, a conformal metric, and a holomorphic cubic differential. It is called
\textnormal{admissible} if
\begin{equation}
K \, = 1 - 2 \, || \Phi ||^2,
\end{equation}
where $|| \Phi ||$ is the norm with respect to $\mathfrak{g}$, and $K$
is the Gaussian curvature of the metric. Let
\begin{equation}
\pi \, : \, \tilde{\Sigma} \, \to \, \Sigma \notag
\end{equation}
be the universal covering of $\, \Sigma$. Then the triple
$( \, \tilde{\Sigma}, \, \pi^* \mathfrak{g}, \, \pi^* \Phi \, )$ is
also admissible and $\tilde{\Sigma}$ admits
an isometric special Legendrian immersion in  $S^5$ with $\pi^* \Phi$
as the associated holomorphic cubic differential.
The immersion is unique up to motion by $SU(3)$.
\end{theorem}
Note that given an admissible triple, there is an $S^{1}$
family of admissible triples by taking $\Phi \to e^{i \tau} \, \Phi$
where $e^{i \tau}$ is any complex number of unit length.

Take a point $p \in \Sigma$  at which $\Phi$ is not 0, and a local
coordinate $z$ in a neighborhood of $\, p\,$ so that
$\Phi = (dz)^3$. A conformal metric
$\mathfrak{g} = e^{2u} dz d\overline{z}$  is compatible, (14), if
\begin{equation}
\Delta u + e^{2u} - 2 \, e^{-4u} = 0.
\end{equation}
(15) is also known as the Monge-Ampere equation describing 2-dimensional
affine spheres. It has the largest, $8$-dimensional, group of symmetry
among nondegenerate(nonlinear) elliptic Monge-Ampere equations.

\begin{definition}
A point $\, p \in \Sigma$ of a special Legendrian surface is
\textnormal{umbilic} if the associated holomorphic differential $\, \Phi$
vanishes at $p$.
\end{definition}

Let $d_p$ denote the degree of zero of $\Phi$ at $p$. Let $S_p$ be
the geodesic $2$-sphere tangent to $\Sigma$ at $p$ and $L_p$ be
the special Lagrangian $3$-plane such that $S_p = L_p \, \cap \, S^5$.
It is clear then that $S_p$ has contact of order $(d_p + 1)$ with
$\Sigma$ at $p$, i.e.,  $(d_p + 1)$th jet of $\phi$ at $p$ is
contained
in $L_p$.

\begin{proposition}
Let $u : \Sigma \to S^5$ be a special Legendrian immersion of
a compact orientable surface of genus g.

1. If $g=0$, $u (\Sigma)$ is the totally geodesic sphere
$S^5 \cap L, \, L \in Slag(3)$.

2. If $g \geq 1$, $\, \sum_{p \in \Sigma} \, d_p = 6g - 6$.
\end{proposition}

\begin{theorem} \textnormal{[Ha]} Let $\Sigma \subset S^5$ be a
special Legendrian surface
with Gaussian curvature $K$.

1. Suppose $K$ is constant. Then either $K = 1$ and $\Sigma$ is
totally geodesic, or $K=0$ and $\Sigma$ is part of the hexagonal torus in Example 1.2.

2. If $\, \Sigma$ is not linearly full in $S^5 \subset \mathbb{R}^6$, it
is totally geodesic.

3. Suppose $\Sigma$ is compact and $K \geq 0$.  Then $\Sigma$ is
either totally geodesic or  the hexagonal torus.
\end{theorem}
\emph{Proof.} $\;$ 1 and 3 follow from $ ||dK||^2 = 16 ||h||^2
||h_1||^2$ and  Simons' type identity
\begin{equation}
-\Delta K = 12K||h||^2 + 8 ||h_1||^2.
\end{equation}
2 is immediate from the structure equations (11), (12), and (13).
$\square$

Special Legendrian surfaces in $S^5$ thus bear resemblance to minimal
surfaces in $S^3$. The extrinsic geometric property, however, can be
significantly different, for special Legendrian
surfaces are integral to the first order differential system (9).

\begin{proposition}
Let $\Sigma \subset S^5$ be a special Legendrian surface,
and $\gamma \subset \Sigma$ a smooth connected curve
($\gamma$ may be a boundary component of $\, \Sigma$).

1. If $\gamma$ contains
an analytic arc, $\, \Sigma$ inherits all the symmetry of  $\, \gamma$.

2. The Gaussian curvature of $\, \Sigma$ along $\gamma$
can be determined from $\, \gamma$ alone.
\end{proposition}
\emph{Proof.}$\;$  It follows from  Proposition 1.2 and the
structure equations (11), (13). $\square$

In terms of the inclusion $\mathbb{C}^3 \subset Im \, \mathbb{O}$,
where  $\mathbb{O}$ is the algebra of Octonians, special Lagrangian
cones are associative cones and their links, special Legendrian surfaces,
are complex curves in $S^6$ with respect to $G_2$ invariant almost
complex structure. Bryant showed every compact Riemann surface admits
a branched immersion in $S^6$ as a complex curve with null torsion
[Br2]. The only special  Legendrian surface with null torsion in
this sense, even locally however, is totally geodesic.

\section{Reflection principle}
Let
\begin{align}
 L^{\perp} &= \{ \, (x, y) \in \mathbb{C}^3 \, | \, x=0 \,\} \, = \,
              \frac{\partial}{\partial y^1} \wedge
\frac{\partial}{\partial y^2} \wedge
              \frac{\partial}{\partial y^3} \notag \\
 \Pi &= \{ \, (z^1, \, z^2, \, z^3 ) \in  \mathbb{C}^3 \, | \, z^3=0
\,\}. \notag
\end{align}
\begin{definition}
By \textnormal{geodesic reflection across $L^{\perp}$} we mean the map
$r_{L^{\perp}} : S^5 \to S^5$ where
\begin{equation}
r_{L^{\perp}} ( x, y) = ( -x, y). \notag
\end{equation}
By \textnormal{geodesic reflection across $\Pi$} we mean the map
$r_{\Pi} : S^5 \to S^5$ where
\begin{equation}
r_{\Pi} ( z^1, \, z^2, \, z^3 ) = ( z^1, \, z^2, \, -z^3 ).\notag
\end{equation}
\end{definition}
Note $r_{L^{\perp}}$ and $r_{\Pi}$ are symmetry of the special
Legendrian differential system (9).

Consider an oriented isotropic 2-plane $\sigma$ in $\mathbb{C}^3$.
Since $SU(3)$ acts transitively on the set of isotropic 2-planes, as well as
on the set of special Lagrangian 3-planes, we may take the 2-plane
\begin{equation}
 \sigma = \, \frac{\partial}{\partial y^1} \wedge
\frac{\partial}{\partial y^2}. \notag
\end{equation}
Then it follows from (9) that
\begin{equation}
L_{\sigma} =  \, - \,  \frac{\partial}{\partial x^3}\ \wedge
        \frac{\partial}{\partial y^1} \wedge \frac{\partial}{\partial
y^2} \notag
\end{equation}
is the unique special Lagrangian 3-plane containing $\sigma$, and the
3-plane
\begin{equation}
L^{\perp} = \,  \frac{\partial}{\partial y^1} \wedge
 \frac{\partial}{\partial y^2} \wedge \frac{\partial}{\partial y^3}
\notag
\end{equation}
is the unique element in $Slag^{\perp}(3)$ that contains $\sigma$.
Let
\begin{equation}
\Pi =  \, \frac{\partial}{\partial x^1} \wedge
\frac{\partial}{\partial y^1} \wedge \frac{\partial}{\partial x^2}
\wedge
\frac{\partial}{\partial y^2} \notag
\end{equation}
be the complex 2-plane generated by $\sigma$. Note  that
both $L^{\perp}$ and $\Pi$ are orthogonal to $L_{\sigma}$ along $\sigma$.

The following reflection principle for special Legendrian surfaces
can be considered as an integral of these linear relations.
\begin{proposition}
1. Let $\Sigma \subset S^{5}$ be a special Legendrian surface
$C^1$ up to its boundary $\, \partial \Sigma$. If part of
$\, \partial \Sigma$ lies in $\,  L^{\perp} \cap S^5 $, then  $\, \Sigma$
is necessarily orthogonal
to $ L^{\perp} \cap S^5$ and can be analytically continued across
each component of $\, \partial \Sigma \cap L^{\perp} $ by
the reflection $r_{L^{\perp}} $. Analogous reflection principle
is true for $\partial \Sigma$ that lies in a complex 2-plane $\Pi$.

2. Let $M \subset \mathbb{C}P^2$ be a Hamiltonian stationary minimizing
Lagrangian surface smooth up to its boundary $\, \partial M$. If part of
$\, \partial M$ lies in $\, \Pi \cap \mathbb{C}P^2$, then
$\, M$ is necessarily orthogonal
to $ \Pi \cap \mathbb{C}P^2$ and can be analytically continued across
each component of $\, \partial M \cap \mathbb{C}P^2 $ by
the reflection $r_{\Pi} $.
\end{proposition}
\emph{Proof.} $\;$ 1. From the orthogonal relation above, the union of
$\Sigma$ and its reflected image is $C^1$ along the corresponding
boundary component. By Proposition 1.3, it is in fact real analytic.

2. $C^1$ Hamiltonian stationary minimizing Lagrangian surface is smooth [ScW].
$\square$

\begin{definition} Let $\Sigma \subset S^{5}$ be a special Legendrian
surface. A curve $\, \gamma \subset \Sigma$ is
\textnormal{anti-special  Lagrangian} if
it lies in $L^{\perp} \in Slag(3)^{\perp}$,
and it is  \textnormal{complex} if it lies in a complex 2-plane.
\end{definition}

An anti-special Lagrangian or  complex curve is necessarily
a geodesic in the analytically continued surface. In fact,
a curve on a special Legendrian surface is anti-special
Lagrangian(complex) if it is a geodesic on which the associated cubic
differential $\Phi$ is real(purely imaginary), Proposition 8.1.

\paragraph{Example 3.1} Consider the hexagonal torus in Example 1.2,
which we parameterize by
$z_j = \frac{1}{\sqrt{3}}  e^{i \, \theta_j}$ with $\theta_1 +
\theta_2 + \theta_3 = 0$.
Each of the curves $\theta_j + \theta_k = const \, $ is anti-special
Lagrangian and each curve $\theta_j - \theta_k = const \,$ is complex.

Reflection principle analogous to Proposition 3.1 can be formulated
for minimal Lagrangian submanifolds in K\"ahler manifolds with
(anti)K\"ahler involution, Section 8.

\section{Rational geodesic polygon}

Let $\Gamma = \{ P_{1}, \, \, .. \, P_{m} = P_{0} \} \subset
S^{5}$, $m \geq 1$,  be a finite union of anti-special Lagrangian
spheres $L^{\perp} \cap S^{5}$ and complex spheres $\Pi \cap S^{5}$.
$\Gamma$ is called a \emph{geodesic m-gon} if $\, \Gamma$ is connected and
its image in $\mathbb{C}P^{2}$ is a $m$-gon, i. e.,
it consists of $m$ vertices and two dimensional faces. We define
$\{ v_{i} = P_{i} \cap P_{i+1} \}_{i=0}^{m-1}$ the vertices of the
polygon with vertex angle $\theta_{i}$ at $v_{i}$, where
$\cos^2{\theta_{i}}$ is the Jacobian of projection from $P_{i}$
to $P_{i+1}$. Note in case the adjacent faces are both complex, the
corresponding vertex is a circle in $S^{5}$.  Let $G_{\Gamma}$ be the
group generated by reflections across each faces of $\Gamma$.

\begin{lemma}
Let $v = P_{1} \cap P_{2} \subset S^{5}$ be a vertex, and
$r_{1}, \, r_{2}$ the reflections across $P_{1}, \, P_{2}$ respectively.
The group generated by $\, \{ \, r_{1}, r_{2} \, \}$ is finite if and only if
the vertex angle at $v$ is a rational multiple of $\, 2 \, \pi$.
\end{lemma}
\noindent \emph{Proof.}$\;$ If $P_{1}, P_{2}$ are both complex or
both anti-special Lagrangian, it follows form the observation that
such pair of planes is \emph{isoclinic} [TeT]. The mixed case
follows similarly, Lemma 4.3. $\square$

\noindent Note the finite group in each case is a Dihedral group.

\begin{definition}
Let $\, \Gamma \subset S^{5}$ be a geodesic polygon. It is a
\emph{rational polygon} if each vertex angle is a rational multiple
of $\; 2 \, \pi$. It is of \emph{finite type} if $\;G_{\Gamma}$ is finite.
\end{definition}
\emph{Remark.} $\;$ Let $X$ be a K\"ahler surface. The above
definition can be extended to a geodesic polygon $\Gamma \subset X$,
where $\Gamma$ is a finite union of the fixed point locus of K\"ahler
or anti-K\"ahler involutions. Note in this case we must require in addition
the two associated angles, [TeT], at every non complex-complex vertex be
rational multiples of $2\, \pi$.

\paragraph{Example 4.1}
Consider a finite group $G \subset U(3)$ generated by unitary
reflections of period two [Cox]. Let $E$ be the union of the axes of
reflections. Any finite ordered set
$\{  \, P_{1}, \, .. \, P_{m}= P_{0} \, \} \subset E$
of distinct elements in $E$ is a rational polygon of finite
type.  A particular class of such polygons will be examined in
detail in Section 6.

We record the following for later application.
\begin{lemma}
Let $\Pi_1$, $\Pi_2$ be complex lines in $\mathbb{C}^2$,
and $L \subset \mathbb{C}^2$ be a Lagrangian plane such that
dim $\Pi \cap L = 1$ for $i=1, \, 2$. Then the angle of the wedge
in $L$ formed by the intersection with $\Pi_1 \cup \Pi_2$ is
determined from the angle between $\Pi_1$, $\Pi_2$.
\end{lemma}
\emph{Proof.} $\;$ Let $\{ \, E_1, \, N_1, \, E_2, \, N_2 \, \}$ be an
orthonormal basis of $\mathbb{C}^2$ with $J(E_i) = N_i$.
Without loss of generality, we may assume $\Pi_1 = E_1 \wedge N_1$, and
$\Pi_2= (\cos(\varphi) E_1 + \sin(\varphi) E_2) \wedge
        (\cos(\varphi) N_1 + \sin(\varphi) N_2)$, where
$\varphi$ is the angle between $\Pi_1$, $\Pi_2$. Suppose $E_1 \in L$.
Then since $L$ is Lagrangian, $L = E_1 \wedge E_2$. $\square$

\begin{lemma}
Let $\,  \gamma_1$, $\, \gamma_2$ be curves on a special Legendrian
surface $\, \Sigma$ that meet at a vertex
$\, p \in \Sigma$ with an angle $\, \varphi$. Let $\, L^{\perp}_1, \,
L^{\perp}_2  \in Slag(3)^{\perp}$ and
$\, \Pi_1, \, \Pi_2  \in  Gr_{\mathbb{C}}(2, \mathbb{C}^3)$.

1. Suppose $\, \gamma_i \subset L^{\perp}_i \cap S^5$ or $\,
\gamma_i \subset \Pi_i \cap S^5$ for $i = 1, 2$.
Then $\, cos^2(\varphi)$ is determined from the angle between $\,
L^{\perp}_1$ and $\, L^{\perp}_2$ or
$\, \Pi_1$ and $\, \Pi_2$.

2. Suppose $\, \gamma_1 \subset L^{\perp}_1 \cap S^5$ and $\,
\gamma_2 \subset \Pi_1 \cap S^5$. Then
$\, sin(2\varphi)$ is determined from the angle between $\,
L^{\perp}_1$ and $\, \Pi_1$.
\end{lemma}
\emph{Proof.} $\;$ We consider the case $\, \gamma_i \subset \Pi_i
\cap S^5$ for $i = 1, \, 2$.
Set $E_i = \frac{\partial}{\partial x^i}$ and  $N_i =
\frac{\partial}{\partial y^i}$.
Let $\Pi_2 = E_2 \wedge N_2 \wedge E_3 \wedge N_3$ and without loss
of generality, let
$E_ 3 = p $ and $E_2$ be the tangent vector of $\gamma_2$ at $p$.
Then $E_1 \wedge E_2$ is the tangent space
of $\Sigma$ at $p$, and tangent vector of $\gamma_1$ at $p$ is of the
form $v = sin(\varphi) E_1 + cos(\varphi) E_2$.
Since $\gamma_1 \subset \Pi_1$, $\Pi_1 = E_3 \wedge N_3 \wedge \, v
\, \wedge J(v)$, and $cos^2(\varphi)$ is the
Jacobian of the orthogonal projection from $\Pi_2$ to $\Pi_1$.
$\square$

\section{Construction procedure}

In this section, we describe a general method toward constructing
complete or compact (branched) minimal Lagrangian surfaces in
compact K\"ahler-Einstein surfaces \emph{without} Lagrangian homotopy classes.
A minimal Lagrangian surface in $\mathbb{C}P^{2}$ then lifts to a
special Legendrian surface in $S^5$, Proposition 1.6. We remark however that main
existence theorem Theorem 5.1 provides only a \emph{Hamiltonian stationary}
Lagrangian disk with at most finitely many
interior singular points, and the disk is not minimal in general.

Let $D$ denote the closed unit disk in the complex plane with the standard
Lebesgue measure, and
$D^{\pm} = \{ (x, y) \in D \, | \, x \gtreqqless 0 \, \}$.
Let $X$ be a compact K\"ahler-Einstein surface without Lagrangian
homotopy classes. Define
\begin{equation}
  W^{1,2}_{L}(D, X) \, =
    \, \{ \, l \in W^{1,2}(D, X) \; | \;
 \; l^{*}(x)\varphi = 0 \;
    \; \mbox{for a. e.  $x \in D$} \; \} \notag
\end{equation}
where $W^{1,2}(D, X)$ is the Sobolev space of maps of
square integrable first derivatives, and $\varphi$ is the symplectic
form on $X$.  For $l \in W^{1,2}_{L}(D, X)$, let $E(l)$ denote
the Dirichlet energy and
$\partial l = l |_{\partial D} \in W^{\frac{1}{2}, 2}(\partial D, X)$.

Let $\Gamma \subset X$ be a rational geodesic polygon.
Define
\begin{align}
 W_{L}^{c} &= \, \{ \, l \in  W^{1,2}_{L}(D, X) \, |
           \, E(l) \leq c \,  \} \notag  \\
 W_{\Gamma} &= \, \{ \, l \in W^{1,2}_{L}(D, X) \, | \,
 \partial l(x) \in \Gamma \; \mbox{for a. e.  $x \in \partial D$} \; \}
 \notag \\
 W_{\Gamma}^{c} &= \, \{ \, l \in  W_{\Gamma} \, |
           \, E(l) \leq c \,  \}. \notag
\end{align}
Then $W^{c}_{L}$ is weakly closed in  $W^{1,2}(D, X)$ [ScW].
For simplicity, $\partial l \subset \Gamma$ would mean
$\partial l(x) \in \Gamma$ for a. e.  $x \in \partial D$.
\begin{lemma}
$W_{\Gamma}^{c}$ is weakly closed in $W^{1,2}_{L}(D, X)$.
\end{lemma}
\noindent \emph{Proof}.$\;$ Suppose $l_{i} \to l$ weakly. By lower semicontinuity
of energy and the arguments in Prop 2.6 [ScW], it suffices to show
$\partial l \ \subset \Gamma$. Since the
trace map $W^{1,2}(D) \to W^{\frac{1}{2},2}(\partial D)$ is
continuous, $\partial l_{i} \to \partial l$ weakly in
$W^{\frac{1}{2},2}(\partial D)$, and a subsequence
$\partial l_{i_{j}} \to \partial l$ in $L^{2}$, and hence a subsequence
$\partial l_{i_{j_{k}}} \to \partial l$ almost everywhere. $\square$

Let $\{ P_{1}, \, \, .. \, P_{m} = P_{0} \}$ be the set of faces
of $\Gamma$. $l \in W_{\Gamma}$ is called \emph{monotone} if
there exists a partition of $\partial D = S^{1}$ by a set of $m$ distinct
points $\{ q_{1}, q_{2}, \, \, .. \, q_{m} = q_{0} \} \subset \partial D$
called \emph{vertices} such that
$\partial l(\overline{q_{i} \, q_{i+1}}) \subset P_{i+1}$ as an $L^{2}$
map for $i=0, \, ... \, m-1$. Each of the arc $\overline{q_{i} \, q_{i+1}}$ is
called an \emph{edge}.

We begin with a relative version of the Gromov-Allcock isoperimetric inequality.
Denote $S^- = \, \{ \, (x, y) \in \partial D^- \, | \, x^2 + y^2 = 1 \} $, and
$I_y =\, \{ \, (x, y) \in \partial D^- \, | \, x = 0 \} $.
\begin{proposition}\textnormal{[Wa2]}
Let $\, \Gamma \subset \mathbb{C}^{2}$ be either a Lagrangian plane,
or a complex plane, or a union of any two of Lagrangian or complex planes
that intersect transversally at the origin.
There exists a constant $\, \mu(\Gamma)$ such that for any $\,
W^{1,2}$ curve $\, C : S^{-} \to \mathbb{C}^{2}$ which bounds a $\, W^{1,2}$
Lagrangian half disk(or a wedge respectively in case $\, \Gamma$
consists of two planes) $\, l_{0} : D^{-} \to  \mathbb{C}^{2}$
with $\, l_{0}(I_y)$ $\subset \Gamma$ and $\, l_{0}|_{S^{-}} =
C$, there exists a $\, W^{1,2}$ Lagrangian half disk(wedge) $\,  l
:  D^{-} \to  \mathbb{C}^{2}$ spanning $\, \Gamma$ and $\, C$ such that
\begin{equation}
    \mbox{Area}(l(D^{-})) \, \leq \mu \, \mbox{Length}^2(C). \notag
\end{equation}
If $\, \Gamma$ contains a complex plane, any curve $\, C$ with boundary in
$\, \Gamma$ always bounds a Lagrangian half disk or a wedge with boundary in
$\, \Gamma$.
\end{proposition}

\noindent These relative isoperimetric inequalities then give
the relative analogues of the Collar Lemma 4.8 in [ScW].

\begin{corollary} Let $\Gamma \subset X$ be a rational
polygon, and suppose $l \in W_{\Gamma}$ is  a weakly conformal,
monotone, weakly Lagrangian map that is minimizing in its relative
homotopy class. Then $l$ is H\"older continuous on $D$ up to boundary.
\end{corollary}
\noindent \emph{Proof}.$\,$ The interior regularity is proved in
Theorem 2.8 in [ScW]. Suppose $q \in \partial D$ is a \emph{vertex},
and consider a small extension $\hat{l}$ of $l$ in a neighborhood of $q$ by
finite successive reflections across the corresponding faces of $\Gamma$,
which is defined since $l$ is monotone and $\Gamma$ is rational. The
corollary follows from the relative isoperimetric inequality and the
interior regularity applied to $\hat{l}$.
Regularity along the edges follows similarly. $\square$

A weakly Lagrangian map $l\in W_{\Gamma}$ is called \emph{stationary} if for any
Lagrangian deformation vector field $V$ along $l$ which is tangent to $\Gamma$ on
$\partial D$, the variation of area
\begin{equation}
    \delta_{V} Area(l(D)) = 0. \notag
\end{equation}
It is \emph{Hamiltonian stationary} if the variation vector field $V$ is assumed to
be Hamiltonian.

\begin{proposition}
Let $\, \Gamma \subset X$ be a rational polygon, and suppose
$\, l \in W_{\Gamma}$
is monotone and stationary. For any $\, q \in \partial D$ consider $\hat{l}$,
which is a slight continuation of $l$ by successive reflections across the
corresponding  face(s) of $\, \Gamma$ at $\, l(q)$. Then
$\, \hat{l}$ is stationary
with respect to compactly supported Lagrangian variations.
\end{proposition}
\noindent \emph{Proof}.
Suppose $q$ is an edge point, and $l(q)$ lies on a face $P$ of
$\Gamma$. Take a representative $u = u^{-} \cup u^{+}$ of $\hat{l}$
near $q$ defined on $D^{-} \cup D^{+}$ so that $q$ corresponds to
the origin and $u^{+} = r_{P} \circ u^{-}$. Given any compactly
supported Lagrangian vector field $V$ along $u$ in a neighborhood
of $q$, let $\tilde{V}$ be the average
of $V$ with respect to $r_{P}$. Then from the symmetry of $u$,
\begin{equation}
   \delta_{V} Area(u(D)) \, = \,
   2 \, \delta_{\tilde{V}} \, Area(u^{-}(D^{-})) \, = 0, \notag
\end{equation}
for $\tilde{V}$ is an admissible variation of $u^{-}$ along
$y-$axis of $D^{-}$. The case $q$ is a vertex follows by a
similar argument. $\square$

We now describe a set of conditions for the rational polygon $\Gamma$ that
is necessary for our construction. For any configuration $Y \subset X$,
$\, \pi^L_2(X, Y) \subset \pi_2(X, Y)$ denotes the set of relative homotopy classes
that can be represented by Lagrangian disks.
\begin{definition}
Let $\, \Gamma \subset X$ be a rational polygon with faces
$\, \{\, P_1, \, ... \, P_m = P_0 \, \}$. $\, \Gamma$ is \emph{admissible} if
the relative Lagrangian homotopy classes $\, \pi_2^L(X, P_i)$,
$\, \pi_2^L(X, P_i \cup P_{i+1}) , \; i = 0, ... \, m-1$, are trivial.
\end{definition}

\paragraph{Example 5.1}  Suppose each faces of a polygon $\, \Gamma \subset X$
is $S^2$, and every adjacent pair of faces intersects at a single point. Then
$\, \Gamma$ is admissible if $X$ is without Lagrangian homotopy class.

Let $\beta \in \pi^L_{2}(X, \Gamma)$ be a nonzero class which can be
represented by a monotone Lagrangian disk. Define
\begin{equation}
    W_{\Gamma, \beta} \, = \, \{ \, l \in W_{\Gamma} \; | \; \,
         [l] = \beta, \quad  l \; \mbox{is monotone}  \}, \notag
\end{equation}
where $[l]$ is the homotopy class represented by $l$. Note the
maps $W_{\Gamma} \to \pi_{2}(X, \Gamma)$ and
$W_{\Gamma} \to W^{\frac{1}{2},2}(\partial D, X) \to \pi_{1}(\Gamma)$
are well defined and continuous [ScU][Ku].

As remarked in [Ku], there exists a lower bound $\delta > 0$ such that
$E(l) \geq \delta$ for any $l \in W_{\Gamma, \beta}$,  $\, \beta \neq 0$.
Set
\begin{align}
    E(\Gamma, \beta) &=  \mbox{inf}_{W_{\Gamma, \beta}}
                         \, E(l)   \notag \\
    A(\Gamma, \beta) &=  \mbox{inf}_{W_{\Gamma, \beta}}
                         \, Area(l). \notag
\end{align}
Then $A(\Gamma, \beta) = \frac{1}{2}  E(\Gamma, \beta)$, and if
$E(l) = E(\Gamma, \beta)$, $l$ is conformal and $A(l) = A(\Gamma, \beta)$.

The following is the main existence theorem of this section.
\begin{theorem} Let $\, \Gamma \subset X$ be a compact admissible complex
rational polygon in a compact K\"ahler-Einstein surface $X$
without Lagrangian homotopy
classes. Let $\beta \in \pi_{2}^L(X, \Gamma)$ be a nonzero class which can be
represented by a monotone Lagrangian disk. There exists an area minimizing
$\, l \in W_{\Gamma, \beta}$  with $A(l) = A(\Gamma, \beta)$ which is a Hamiltonian
stationary branched Lagrangian surface with at most finitely many interior conical
singular points and smooth up to $\partial D$. $\, l(D)$ generates a complete Lagrangian
surface by successive reflections across its boundary. If $\, l(D)$ is smooth,
it is minimal.
\end{theorem}
\noindent \emph{Proof}. $\;$  Let $\{l_{k}\}_{k=1}^{\infty} \subset W_{\Gamma, \beta}$ be
a minimizing sequence of monotone Lagrangian maps with vertices
$\{ \, q_k^j \}_{j=1}^m$ that weakly converges to $l \in W_{\Gamma}$.

\noindent Case $m \leq 3$ or $\Gamma$ is at most a 3-gon : $\;$
In this case, we may assume each $q_k^j = q^j$, $j= 1, \, ... \, m$, is
fixed for all $k$, and hence $l$ is also monotone, proof of Lemma 5.1.

Interior regularity. $\;$ It follows from Theorem 4.10, Proposition 5.7 of [ScW].

Boundary regularity.

1. Since $\pi_2^L(X, P_j) = 0$ and $\pi_2^L(X, P_{j} \cup P_{j+1}) = 0$, we can
apply the relative collar lemma and reflection to show there exist at most finitely
many edge points and vertices at which the energy measure of the sequence
$\,\{ \, l_k \, \}$
blow up. Each blow up points can then be used to construct a sequence of minimizing
Lagrangian disks(edge point) or wedges(vertex) by the relative isoperimetric
inequality. Again since $\Gamma$ is admissible, such sequence is trivial, and the
blow up loci in $\partial D$ is empty and $l_k \to l$ strongly in $W^{1,2}$.

2. Let $q \in \partial D$ be an edge point(or a vertex respectively), and let
$\hat{l}$ be the continuation of $l$  in a neighborhood of $q$ via
(successive) reflection(s) across the face(s) of $\Gamma$
that contain $l(q)$. The relative isoperimetric inequality
and collar lemma gives a relative version of the strong
$W^{1,2}_{loc}$ compactness theorem Proposition 4.7 of [ScW], which then shows
the existence of the tangent cone of $\hat{l}$ at $q$.

3. From 1 and 2 and reflectional symmetry, the link of each tangent cone of
$\hat{l}$ along the edge points
of $\partial D$ is $C^1$, and from the normal form of the Hamiltonian stationary
Lagrangian cones in $\mathbb{C}^2$, [ScW], it is in fact smooth.
One can easily check that the only
Hamiltonian stationary Lagrangian cone with reflectional symmetry across a
complex plane is flat, and hence
each tangent map is proper holomorphic with respect to a suitable coordinates
via Schwartz reflection principle.
Now boundary regularity follows as in Theorem 4.10 of [ScW] using partial
regularity Theorem 4.1 of [ScW]. Regularity at the vertices
is verified similarly.

\noindent Case $m \geq 4$ : $\;$
From the arguments in Case $\, m \leq 3$, it suffices to show $l$ is monotone.

Suppose $m=4$, and $q_k^j = q^j$ is fixed for $j =$ 1, 2, 3. If the sequence
$q_k^4 \to q^1$, we claim that the energy measure of the sequence $l_k$ blows up at
$q^1$.  For if not, $l_k \to l$ strongly in $W^{1,2}$ in a neighborhood of $q^1$,
and there exists a neighborhood $U$ of $q^1$ for any neighborhood $V$ of $l(q^1)$ in $X$
such that $l_k(U) \subset V$ for all sufficiently large $k$. Since each $l_k$ is
monotone, this is a contradiction if we choose $V$ so small that
$V \cap P_0 = \emptyset$.

We now apply arguments in 1. of Case $\, m \leq 3$ above to show there exists a
minimizing sequence of Lagrangian wedges spanning $P_0 \cup P_1$,
which is trivial for $\pi^L_2(X, P_0 \cup P_1)$ is trivial. Hence the sequence
$\{ \, q_k^4 \,\}$ has a limit point in $\overline{q^3 \, q^1} \subset \partial D$
distinct from $q^1$ and $q^3$, and thus $l$ is monotone.

The case $m >4$ can be proved by successive application of the
arguments for the case $m=4$.

If $l$ is an immersion up to boundary, the mean curvature deformation
is admissible along the boundary due to reflectional symmetry, and
$l$ is minimal. General case with branch points follows
similarly by applying Lemma 8.2 in [ScW].  $\square$

\noindent \emph{Remark.} $\, $ For general K\"ahler-Einstein surface $\, X$ and
a compact rational complex polygon $\, \Gamma$, a relative Lagrangian
homotopy or homology class with respect to $\, \Gamma$ can be represented
by a finite union of  Hamiltonian stationary Lagrangian surfaces,
wedges and disks with boundary in $\, \Gamma$,  and  Lagrangian spheres.

\paragraph{Example 5.2}
Let $\, \Gamma \subset \mathbb{C}^2$ be a rational complex polygon,
and let $\alpha_{\Gamma} \in \pi_1({\Gamma})$ be the class represented by
a monotone closed curve that connects the vertices of $\, \Gamma$.
By considering a sufficiently large torus
$\mathbb{C}^2/\mathbb{Z}$ for some lattice $\mathbb{Z}$,  Proposition 5.1
and Theorem 5.1 imply there exists a monotone Hamiltonian stationary
Lagrangian disk $\, l_{\Gamma}$ with $\, \partial l_{\Gamma} \subset \Gamma$
such that $\, [\partial l_{\Gamma}] = \alpha_{\Gamma}$.

Suppose $\, l_{\Gamma}$ is smooth, hence minimal.
Let $\, dz^1 \wedge dz^2$ be the holomorphic volume form
of $\mathbb{C}^2$, and, with a slight abuse of notation, consider
\begin{equation}
A(\Gamma)\, e^{-i \, \Theta(\Gamma)} = \frac{1}{2} \,
\int_{\alpha_{\Gamma}} \, z^1 \, dz^2 - \, z^2 \,dz^1. \notag
\end{equation}
Since $\, dz^1 \wedge dz^2 = 0$  on $\, \Gamma$, the above integral
is well defined. $\, A(\Gamma) $ and $\, \Theta(\Gamma)$ represent
the area and the \emph{phase} of the minimal Lagrangian
surface $\, l_{\Gamma}$ respectively.

\begin{corollary}
Let $X$ be a compact K\"ahler-Einstein surface. Let $\, \sigma$ be a
K\"ahler involution, and $\, \Gamma_{\sigma}$ be a compact
fixed point locus of $\, \sigma$, which is a complex curve.
Suppose $\, \pi^L_2(X, \Gamma_{\sigma})$ is nontrivial.
Then there exists a Hamiltonian stationary Lagrangian sphere in $\, X$,
which is smooth except at most finitely many conical singularities.
\end{corollary}
\noindent \emph{Proof.} \emph{Remark} below Theorem 5.1. $\square$

\begin{corollary}
Let $\, \Gamma \subset \mathbb{C}P^{2}$ be a rational complex polygon,
and let $\, \alpha \in \pi_1(\Gamma)$
be the class which is represented by a monotone curve that connects the vertices
of $\, \Gamma$. Then there exists a Hamiltonian stationary Lagrangian disk
$l \in W_{\Gamma}$ with $[\partial l] = \alpha \in \pi_{1}(\Gamma)$ which
is smooth up to boundary except at most finitely many interior
conical singularities.
\end{corollary}
\noindent \emph{Proof.}
Let $\, \gamma_0 \in \Gamma$ be a monotone curve with $\, [\gamma_0] = \alpha$.
By local accessibility theorem for contact structures, [BCG], there exists
a Legendrian lift $\, \gamma \subset S^5$ of $\, \gamma_0$.
Let $\, D_{0} \subset S^{5}$ be a disk with boundary $\, \gamma$.
Then by triangulization and local accessibility theorem for contact structures
again, we can deform $\, D_{0}$ to a Legendrian disk while keeping the boundary
fixed. Since $\, \Gamma$ is admissible by Example 5.1, the
corollary follows from Theorem 5.1 and the fact the trace map
$W^{1,2}(D, \Gamma) \to W^{\frac{1}{2},2}(\partial D, \Gamma)$ is
continuous. $\square$

Let $\, M^0_{\Gamma, \alpha}$ be the Hamiltonian stationary disk in the above corollary.
Then by successive reflection across each edges, we obtain a complete
Hamiltonian stationary Lagrangian surface $\, M_{\Gamma, \alpha}$ in
$\, \mathbb{C}P^2$, or a compact surface if $\, G_{\Gamma}$ is finite. Let
$\, \Sigma_{\Gamma, \alpha} \subset S^5$ be a connected  Legendrian lift of
$\, M_{\Gamma, \alpha}$. $\, \Sigma_{\Gamma, \alpha} \to M_{\Gamma, \alpha}$ is a
nontrivial covering when $\, G_{\Gamma} \cap Z(U(3)) \ne I_3$, where
$Z(U(3))$ is the center of $U(3)$. For instance,  $M_{\Gamma, \beta}$ is
nonorientable whenever $-I_3 \in G_{\Gamma}$, Corollary 1.1.

\section{ Surfaces $\, \Sigma_{k, 3}$}

Let $k \geq 3$ be an integer and set
\begin{equation}
    \epsilon = \exp{\frac{2\pi i}{k}}. \notag
\end{equation}
Consider a set of lines in $\mathbb{C}P^2$
\begin{align}
    P_{1} &= \{ \, [z_{1}, \, z_{2}, \,z_{3} ]\in \mathbb{C}P^2 \, | \;
                  z_{1} = \epsilon \, z_{2}  \, \} \notag \\
    P_{2} &= \{ \, [z_{1}, \, z_{2}, \,z_{3} ]\in \mathbb{C}P^2 \, | \;
                  z_{1} = z_{2} \, \} \notag \\
    P_{3} &= \{ \, [z_{1}, \, z_{2}, \,z_{3} ]\in \mathbb{C}P^2 \, | \;
                  z_{2} =  z_{3} \, \}, \notag
\end{align}
and let $\{r_{i}\}_{i=1}^3$ be the corresponding reflections. Let
$\, \Gamma = \cup_{i=1}^3 P_i$ be the complex geodesic triangle with vertices
\begin{align}
    v_{i} &= P_{j} \cap  P_{k} \; \; (ijk)=(123), \notag
\end{align}
and the set of vertex angles
\begin{equation}
    \{ \; \frac{\pi}{3}, \, \frac{\pi}{3}, \, \frac{\pi}{k} \, \}.
    \notag
\end{equation}
The group $G_{k, 3}$ generated by reflections $\{r_{i}\}$ is finite of order
$6k^{2}$ with defining relations
\begin{align}
    r_{i}^{2} &= (r_{2}r_{3})^{3}= (r_{3}r_{1})^{3}=
            (r_{1}r_{2})^{k}=(r_{1}r_{2}r_{3}r_{2})^{3}= 1. \notag
\end{align}
The image of  $\, \Gamma$ under successive reflections across
edges gives rise to a tessellation of $\mathbb{C}P^{2}$ with $6k^{2}$
triangular faces and $3k+2k^{2}$ vertices [Cox].

Let $\, M_{\Gamma}$ be a minimizer of area among all Lagrangian
disks spanning $\, \Gamma$, and let $\, M_{k, 3}$ be the compact Hamiltonian
stationary Lagrangian surface obtained by successive reflection of
$\, M_{\Gamma}$. If $\, M_{\Gamma}$ is an immersed disk, we get from
Gauss-Bonnet formula(since $\, M_{\Gamma}$ spans the smaller angles at
each vertex by the minimizing property) and Lemma 4.2,
\begin{equation}
    \int_{M_{k, 3}} K \, dA = \, 6 \, k^{2} \,
           \pi \, (\frac{1}{k} - \frac{1}{3}), \notag
\end{equation}
and hence
\begin{equation}
    \chi(M_{k, 3}) = \, k(3-k). \notag
\end{equation}

Consider now the Legendrian lift $\, \Sigma_{k, 3} \subset S^{5}$. Since the center
of $\, G_{k, 3}$ is cyclic of order $\, (3, k)$, the greatest common divisor of
3 and $k$,
\begin{equation}
\Sigma_{k, 3} \to M_{k, 3} \notag
\end{equation}
is 3:1 when $k$ is a multiple of 3, and otherwise 1:1.
\begin{theorem}
For each integer $k \geq 3$, there exist a compact Hamiltonian stationary
Legendrian surface
$ \, \Sigma_{k, 3} \subset S^5$ of genus $\, 3(1+\frac{k(k-3)}{2})$ if $3 \, | \, k$,
and genus $1+\frac{k(k-3)}{2}$ otherwise. $ \, \Sigma_{k, 3}$ is smooth
except at most finitely many conical singularities.
The image of this surface under Hopf map
is a compact orientable Hamiltonian stationary Lagrangian surface in
$\, \mathbb{C}P^2$ of genus $1+\frac{k(k-3)}{2}$.
\end{theorem}
\noindent \emph{Remark.} Compact orientable surface of genus $\geq 2$
cannot be embedded in $\, \mathbb{C}P^2$ as a Lagrangian surface by
Proposition 1.8.

$G_{k, 3}$ acts on the vertices of the tessellation mentioned
above with three orbits $V_{+}$, $V_{-}$, and $V_{0}$ of
order $k^{2}$, $k^{2}$, and $3k$ respectively. Moreover, there exists a Lagrangian
reflection of $\mathbb{C}P^{2}$ that interchanges  $V_{+}$, $V_{-}$,
and leaves $V_{0}$ invariant as a set [Cox]. Note this corresponds to adding
a single generator of order two to $G_{k, 3}$. The enlarged group is also finite,
and hence we obtain a rational Lagrangian polyhedral variety that consists of
equilateral triangles of vertex angle $\frac{\pi}{k}$. From the symmetry, it is
likely that $\Sigma_{k, 3}$ is invariant under these enlarged
group and hence can be constructed by successive reflections of a fundamental
domain spanning the equilateral Lagrangian triangle.

Suppose $M_{k, 3}$ is smooth, hence minimal, and consider the subgroup
$G_{k, 3}^{even} < G_{k, 3}$ of even elements, which is well
defined for the defining relations are all even. From Hurwitz formula,
\begin{equation}
   M_{k, 3} \to  M_{k, 3}/G_{k, 3}^{even} = \mathbb{C}P^{1} \notag
\end{equation}
is a $3k^{2}$-fold cover branched over $V_{+}$, $V_{-}$, and $V_{0}$
with branching degrees 3, 3, and $k$ respectively. This leads to an
alternative description of $M_{k,3}$ in terms of a singular solution to
(15) on a hexagonal torus as follows. The original idea of quotient by
discrete symmetry is due to Robert Bryant in his unpublished note on
Lawson's minimal surfaces. We follow the notations adopted in Section 2.

Let $w$ be a coordinate of $M_{k, 3}/G_{k, 3}^{even}$, and set
\begin{align}
    w(V_{\pm}) &= \, \pm 1 \notag \\
    w(V_{0}) &= \, 0. \notag
\end{align}
As $\,\{ \,  V_{+}, V_{-}, V_{0} \, \} $ is the umbilic locus of $M_{k,3}$, the holomorphic
cubic differential $\Phi$ of $M_{k,3}$, Section 2, is then a constant multiple of
\begin{equation}
    \frac{(dw)^{3}}{w^{2}(w^{2}-1)^{2}}. \notag
\end{equation}

Since $G_{k, 3}^{even}$ acts as isometry, holomorphic cubic differential $\Phi$ and
the metric $\mathfrak{g}$ on $M_{k, 3}$
can be pushed forward to $M_{k, 3}/G_{k, 3}^{even}$ as
a meromorphic differential $\Phi_{0}$ and
a singular metric $\mathfrak{g}_{0}$.  Let $T$ be the
hexagonal torus in Example 1.2, and
$\pi_3 :  T \to \mathbb{C}P^1$ be the threefold covering branched at
three points $\, \{  \, \pi_3^{-1}(1),  \, \pi_3^{-1}(0),
 \, \pi_3^{-1}(1) \, \}$.
Let $z$ be the standard coordinate on the universal covering
$\mathbb{C} \to  T$.  Then $\pi_3^* (\Phi_{0})$ is a holomorphic
differential on $T$. In fact, we may choose $z$ so that
\begin{equation}
\pi_3^* (\Phi_0) = \lambda (dz)^3  \notag
\end{equation}
for some real $\lambda > 0$. Similar analysis for the metric
$\mathfrak{g}_0$ shows if we write
\begin{equation}
\pi_3^* (\mathfrak{g}_0) = e^{2u} dz \, d\tilde{z}, \notag
\end{equation}
then $u(z, \overline{z})$ is a function on $T$
that  satisfies the compatibility condition
\begin{equation}
\triangle \, u \, + e^{2u} - 2 \lambda^2 e^{-4u} \, = \, 0 \notag
\end{equation}
with appropriate logarithmic singularities at the branch points of
$\pi_3$.

\section{Gauss map, polar surface, and bipolar surface}
We now turn our attention to surfaces associated to a given special
Legendrian surface.
It will be shown for example that generating functions of the contact
transformations of $S^5$ induced by $SU(3)$ action,
when restricted to a special Legendrian surface, describe a minimal
surface in the unit sphere
of the Lie algebra $\mathfrak{su}(3)$. To draw geometric conclusion
from this, we apply
maximum principle to obtain the following  half space theorem for
compact minimal Lagrangian surfaces in $\mathbb{C}P^2$.

\emph{ A compact minimal Lagrangian surface in $\mathbb{C}P^2$ does
not lie in any open geodesic ball of
radius $\frac{\pi}{2} -  arccos(\frac{1}{\sqrt{3}})$. }

We continue to use the notation adopted in Section 2.

\subsection{Gauss map}

Let $u : \Sigma \to S^5 $ be a special Legendrian surface.
Gauss map $u^* : \Sigma \to Isot^{+}(2, \mathbb{C}^3) =
SU(3)/SO(2)$
 is defined by
\begin{equation}
u^* = e_1 \wedge e_2.
\end{equation}
$u^*$ is then a conformal and minimal immersion of $\Sigma$ with the
induced metric
\begin{equation}
d u^* \circ d u^* = ( 2 - K ) d u \circ d u. \notag
\end{equation}

$Isot^{+}(2, \mathbb{C}^3)$ admits an integrable  $CR$-structure as a
real hypersurface of the hyperquadric
$Gr^+_{\mathbb{R}}(2, \mathbb{C}^3) \subset \mathbb{C}P^5$.
Gauss map $u^*$ is never a complex curve with respect to this
$CR$-structure unless
it is totally geodesic. Instead, consider the fibration in Lemma 1.1
\begin{equation}
Isot^{+}(2, \mathbb{C}^3) \to S^5  \notag
\end{equation}
with $\mathbb{C}P^1$ fibers. Upon reversing the orientation of each
fiber, we obtain a new
$CR$-structure on $Isot^{+}(2, \mathbb{C}^3)$, under which the Gauss
map is complex.
This $\, CR$-structure is however not integrable [Sal].


Let $\gamma : S^1 \to S^5$ be an isotropic curve.
Since the differential system for special Legendrian surfaces (9)
is invariant under $SU(3)$ action,
Noether's theorem implies  not every such $\gamma$ bounds a special
Legendrian surface(variety or current).
Equivalently, the first order characteristic cohomology of the
special Legendrian differential system is at
least of dimension $8$.  Now by Lemma 1.1, every isotropic curve
$\gamma$ has a unique lift
$\gamma^* : S^1 \to Isot^{+}(2, \mathbb{C}^3)$. Thus, a necessary
condition for $\gamma$ to bound
a special Legendrian surface is that its lift $\gamma^*$  must bound
a complex variety in $Isot^{+}(2, \mathbb{C}^3)$.
More specifically, $\gamma^*$ must satisfy
\begin{equation}
\int_{\gamma^*} \eta = 0, \notag
\end{equation}
where $\eta$ is a 1-form on $Isot^{+}(2, \mathbb{C}^3)$ such that $d
\eta $ modulo contact form is a linear
combination of forms of type $(2,0)$ or $(0,2)$ with respect to the
aforementioned nonintegrable $CR$-structure.

In case of a curve in a complex manifold, a necessary and sufficient
condition for a curve to
bound a complex variety is known, and has been generalized to higher
dimensions
by Harvey and Lawson [HL2].

\subsection{Polar surface}

We now view (17) as a map $u^* : \Sigma \to \bigwedge^2
\mathbb{R}^6$.
Consider the endomorphism of $\bigwedge^2 \mathbb{R}^6$ induced from
the complex structure $J$ on $\mathbb{C}^3 = \mathbb{R}^6$,
which we continue to denote by $J$. Then
$\, J \circ J = 1_{\bigwedge^2 \mathbb{R}^6}$, and let
$\, \bigwedge^2 \mathbb{R}^6 = \bigwedge^9_{+} \oplus
\bigwedge^6_{-}$ be the $\pm 1$ eigenspace
decomposition. $\bigwedge^9_{+}$ further decomposes as
$\bigwedge^9_{+} = \mathbb{R} \oplus W^8$ where $W$ is the adjoint
representation of $SU(3)$.

Let
\begin{align}
u^*_+ &= \frac{1}{\sqrt{2}}( e_1 \wedge e_2 + n_1 \wedge n_2 )
\subset S^7 \subset W^8  \notag \\
u^*_- &= \frac{1}{\sqrt{2}}( e_1 \wedge e_2 - n_1 \wedge n_2 )
\subset S^5 \subset \mbox{$\bigwedge^6_{-}$} \notag
\end{align}
be the respective projections, and note that $u^*_{-}$ is
identified with $u$ under the isomorphism
$\mathbb{C}^3 = \bigwedge^2 \mathbb{C}^3$. From the structure
equations in Section 2,
\begin{align}
d u^*_{+} \circ  d u^*_{+} &= ( 3 - 2 K  ) \, d u \circ d u
 \notag  \\
\Delta \, u^*_{+} &= - 2 \, ( 3 - 2 K ) \, u^*_{+}    \\
\Delta \, u^*_{-} &= - 2 \, u^*_{-}.
\end{align}
Thus $u^*_{+}$ is a minimal immersion into $S^7$, and we call
$u^*_{+}$
\emph{polar immersion} of the special Legendrian immersion $u$.
Note that $u^*_{+}$ is well defined on the minimal Lagrangian
image of $\Sigma$ under Hopf map.

Let $\xi \in \bigwedge^2 \mathbb{R}^6$ be a 2-vector, and set
\begin{equation}
u^{\xi}_{\pm} = \, \sqrt{2} \langle  u^*_{\pm}, \, \xi \rangle \,
= \, \langle e_1 \wedge e_2 ,  \, \xi \pm J(\xi) \rangle. \notag
\end{equation}
From (18) and (19), we get
\begin{align}
\Delta \, u^{\xi}_+ &= -2 \, (3 - 2 K ) \, u^{\xi}_+  \notag \\
\Delta \, u^{\xi}_- &= -2 \, u^{\xi}_- . \notag
\end{align}

\begin{proposition}
Let $\, \Sigma \subset S^{5}$ be a compact special Legendrian surface,
and let $\, e_1 \wedge e_2$ be its oriented unit tangent plane
field. For any 2-vector $\, \xi \in \bigwedge^2 \mathbb{R}^6$,
\begin{align}
\int _{\Sigma} ( 3- 2 K ) \, \langle e_1 \wedge e_2 ,  \, \xi +
J(\xi) \rangle \, dvol_{\Sigma} &= 0  \notag \\
\int _{\Sigma}  \langle e_1 \wedge e_2 ,  \, \xi - J(\xi) \rangle  \,
dvol_{\Sigma} &= 0. \notag
\end{align}
In particular, for any  2-vector $\, \xi \in W^8$ or $\,  \xi \in
\bigwedge^6_-$, there exists a point $\, p \in \Sigma$
such that $\, \langle e_1 \wedge e_2 |_p, \, \xi \rangle = 0$.
\end{proposition}

It is easily verified that if $\Sigma$ is compact with genus $\geq
2$, then its polar surface is
linearly full in $W^8$.

\subsection{Bipolar surface}

Let $\varpi^{\sharp} = e_ 1\wedge n_1 +  e_2 \wedge n_2 + e_3 \wedge
n_3$ be the 2-vector dual to the K\"ahler  form $\varpi$ of $\mathbb{C}
^3$, which is constant along $\Sigma$. We define
the \emph{bipolar immersion} associated to a special Legendrian
immersion $u : \Sigma \to S^5$ by
\begin{equation}
u_* = \sqrt{\frac{3}{2}} \, ( e_3 \wedge n_3 - \frac{1}{3}
\varpi^{\sharp} ) \subset S^7 \subset W^8.
\end{equation}
It is a conformal immersion with
\begin{align}
d u_* \circ d u_* &= 3 \, du \circ d u  \notag \\
\Delta \, u_* &= - 6 \, u_*.
\end{align}
Thus the bipolar immersion associated to a special Legendrian
immersion is also minimal.
$u_*$ is also well defined on the minimal Lagrangian image of
$\Sigma$ under the Hopf map.

Coordinate functions of bipolar immersion admit a simple geometrical
interpretation in terms of
generating functions.
Let $\xi \in \mathfrak{su}(3) = W^8$, and denote its induced vector
field on $S^5$ by $V_{\xi}$.
Define the generating functions $G_{\xi} \in C^{\infty}(S^5)$ by
\begin{align}
 G_{\xi}(p) &=  \langle p \wedge  J(p) , \, \xi \rangle    \notag  \\
            &=  \langle V_{\xi}(p), \,   J(p)  \rangle. \notag
\end{align}
Each $G_{\xi}$ is again a well defined function on $\mathbb{C}P^2$.
\begin{lemma}
$\, G_{\xi} : \mathfrak{su}(3)  \to  C^{\infty}(S^5)$ is one to one.
\end{lemma}
\emph{Proof.} $\;$ Suppose $\langle V_{\xi}(p), \,   J(p)  \rangle
\equiv 0$.  Let $\theta$ be the contact form (2), and
$\mathfrak{L}_{V_{\xi}} \theta$ the Lie derivative of $\theta$.
The lemma follows from the equation
\begin{equation}
\mathfrak{L}_{V_{\xi}} \theta = 0 = d(V_{\xi} \, \lrcorner \, \theta)
                     + V_{\xi} \, \lrcorner \, d\theta =  V_{\xi} \,
\lrcorner \, d\theta,
\end{equation}
for $d\theta$ is nondegenerate. $\square$

Let $\Sigma$ be a special Legendrian surface in $S^5$. Put
\begin{align}
u_{\xi} &= \, \sqrt{\frac{2}{3}} \langle \, u_*,  \, \xi \,
\rangle \notag  \\
        &= \; G_{\xi} \, \circ \, u  \notag
\end{align}
and denote $\mathfrak{G}(\Sigma) = \{ \, u_{\xi} \; | \; \, \xi \in
\mathfrak{su}(3) \, \}$.
Then
\begin{equation}
    \Delta u_{\xi} = -6 \, u_{\xi}, \notag
\end{equation}
by (21).

Lemma 7,1 implies $\mathfrak{G}(\Sigma)$ can be identified with the
set of Killing - Jacobi fields on $\Sigma$, i.e., the Jacobi
fields along $\Sigma$ obtained by orthogonal projection to the normal
bundle $N\Sigma$ of Killing fields generated by $SU(3)$ action.
\begin{proposition}
Let $\, \xi \in \mathfrak{su}(3)$. Then  $\,  u_{\xi} \equiv 0 \; \,
\mbox{on} \; \Sigma$ if and only if $\, V_{\xi}$ is tangent to $\,
\Sigma$.
\end{proposition}
\emph{Proof.} $\;$ From the structure equations of the homogeneous
space $S^5 = SU(3)/SU(2)$,
there exists a set of 1-forms  $\omega^1, \omega^2, \theta_1, \theta_2$
in a neighborhood of a point $p \in \Sigma$  so that
$-d\theta = \omega^1 \wedge \theta^1 + \omega^2 \wedge \theta^2$ and
$\theta_1, \theta_2 = 0$ on $\Sigma$.
The proposition follows from (22). $\square$

\begin{theorem}
Let $\, \Sigma$ be a  special Legendrian surface in $S^5$. Then dim
$\, \mathfrak{G}(\Sigma) \geq 5$,
and either

dim $\, \mathfrak{G}(\Sigma) = 5$ and $\, \Sigma$ is totally geodesic, or

dim $\, \mathfrak{G}(\Sigma) = 6$ and $\, \Sigma$ is a part of the
hexagonal torus, or

dim $\, \mathfrak{G}(\Sigma) = 7$ or $\, 8$.

\noindent If $\, \Sigma$ is compact and dim $ \mathfrak{G}(\Sigma) =
7$, $\, \Sigma$ is a torus. If $\, \Sigma$ is compact
with genus $\, \geq 2$, then dim $\mathfrak{G}(\Sigma) = 8$.
\end{theorem}
\begin{corollary}
Let $\Sigma$ be a compact special Legendrian surface of genus $\,
\geq 2$. Then the eigenspace of $\, -\Delta$ on $\, C^{\infty}(\Sigma)$
with eigenvalue $6$ is of dimension at least $8$.
\end{corollary}

$\mathfrak{G}(M)$ can be defined for a compact
minimal Lagrangian submanifold $M$  in $\mathbb{C}P^n$.
For instance, dim $\mathfrak{G}(M) \geq n(n+3)/2$ with equality only
if $M$ is totally geodesic. $\mathfrak{G}(M)$
is a subspace of the eigenspace of $- \Delta$ with eigenvalue $2n+2$.

Equation (21) also has the following geometric consequences when
integrated over a compact special Legendrian
surface. Note  (21) can be rewritten as
\begin{equation}
\Delta \, (e_3 \wedge J(e_3)) = - 6 \, e_3 \wedge J(e_3) + 2
\varpi^{\sharp}. \notag
\end{equation}

\begin{lemma}
Let $\, v', v$ be vectors in  $\, \mathbb{R}^{2n} = \mathbb{C}^n$
with complex structure $J$. Then
\begin{align}
\langle v' \wedge J(v'), \, v \wedge J(v) \rangle &= \langle v, \,
v' \rangle^2 \,  + \, \langle Jv, \, v' \rangle^2  \notag \\
                                                  &= \langle  v, \,
v' \rangle^2  \, + \, \langle v, \, J(v') \rangle^2. \notag
\end{align}
\end{lemma}

\begin{proposition}
Let $\,  e_3 : \Sigma \to S^5 \subset \mathbb{C}^3$ be a
special Legendrian immersion  of a compact surface $\, \Sigma$. For
any
vector $\, v \in \mathbb{C}^3$ of unit length, set
\begin{equation}
u_v = \, \langle e_3 \wedge J(e_3),  \, v \wedge J(v) \rangle \, =
                \, \langle e_3, \, v \rangle^2 \, + \, \langle
J(e_3), \, v \rangle^2, \notag
\end{equation}
which represents the length square of the orthogonal projection of
$\, v$ onto the  2-plane $\, e_3 \wedge J(e_3)$.
Then
\begin{equation}
\Delta \, u_v = -6 \, u_v + 2. \notag
\end{equation}
Upon integration,
\begin{equation}
\int_{\Sigma} \, ( u_v - \frac{1}{3} ) \,  dvol_{\Sigma} = 0. \notag
\end{equation}
\end{proposition}

\begin{theorem}
Let $M \subset \mathbb{C}P^2$ be a compact minimal Lagrangian surface,
and let $dist(q, M)$ denote the distance of a point $q \in
\mathbb{C}P^2$ to $M$.
Then
\begin{equation}
dist(q, M) \leq arccos( \frac{1}{\sqrt{3}} ) \notag
\end{equation}
for any $\, q \in \mathbb{C}P^2$.
\end{theorem}
\emph{Proof.} $\;$ Take a point  $\, \tilde{q} \in S^5$ in the inverse
image of $\, q \in \mathbb{C}P^2$ under Hopf map.
By Proposition 7.3, there exists a point $\, p \in M $ such that
$\, dist( \tilde{q}, \pi^{-1}(p) ) = arccos( \frac{1}{\sqrt{3}} )$.
The theorem follows for Hopf map is an isometric
submersion, and hence distance nonincreasing. $\square$

\begin{corollary}
Let $M \subset \mathbb{C}P^2$ be a compact minimal Lagrangian
surface. Then
$M$ does not lie in any open  geodesic ball of radius $\,
\frac{\pi}{2} - arccos (\frac{1}{\sqrt{3}})$.
\end{corollary}
\emph{Proof.} $\;$ Let $ q, \, q' \in \mathbb{C}P^2$ be such that $\,
dist(q, q') = \frac{\pi}{2}$.
The open geodesic ball of radius $\, \frac{\pi}{2} - arccos
(\frac{1}{\sqrt{3}})$ centered at $\, q\,$ is disjoint
from the closed geodesic ball of radius $\, arccos
(\frac{1}{\sqrt{3}})$ centered at $\, q'$. $\square$

Theorem 7.2 agrees with Corollary 1.2  in that the area of a compact
minimal Lagrangian surface cannot be too small.
We also mention Theorem 7.2 and Corollary 7.2 is true for compact
minimal Lagrangian submanifolds in
$\mathbb{C}P^n, \, n \geq 2$, with the constant
$arccos(\frac{1}{\sqrt{n+1}})$.


For minimal Lagrangian submanifolds in K\"ahler-Einstein manifolds,
there is a notion of
Hamiltonian stability and index introduced by Oh, which is finer or
weaker than the usual stability and index for minimal
submanifolds [Oh]. Since submanifolds in
Example 1.1 are all Hamiltonian stable, [AO],
we have many Hamiltonian stable compact  minimal
Lagrangian submanifolds in $\mathbb{C}P^n$ via product, Proposition 1.5,
whereas a compact minimal submanifold in $\mathbb{C}P^n$ other than
complex varieties is unstable [LS].

A compact minimal Lagrangian submanifold $M$ in  $\mathbb{C}P^n$ is
stable under Hamiltonian deformations near  identity
if $\lambda_1 \geq 2(n+1)$, where $\lambda_1$ is the first nonzero
eigenvalue of $\, -\Delta$ on $C^{\infty}(M)$ [Oh].
The remark below Corollary 7.1 gives a slight strengthening of
this result.
\begin{theorem}
A compact minimal Lagrangian submanifold $\,M$ in $\, \mathbb{C}P^n$
is (locally) Hamiltonian stable if
$ \, \lambda_1 = 2(n+1)$, where $\,\lambda_1$ is the first nonzero
eigenvalue  of $\, -\Delta$ on $\, C^{\infty}(M)$.
\end{theorem}

\noindent \emph{Remark}. Are examples (4), (5), (6), (7), (8) and their
product the only Hamiltonian stable compact minimal Lagrangian
submanifolds in  complex projective spaces ?

\section{Intrinsic characterization and dual reflection\\ principle}

In this section,  we consider minimal Lagrangian surfaces in a
complex space form $X(c)$
of constant holomorphic curvature $4c$, $c= 1, 0$ or $-1$ for unified
treatments. Special Legendrian surfaces in $S^5$
are naturally identified with minimal Lagrangian surfaces in $X(1)$.

We first recall the fundamental theorem for minimal Lagrangian
surfaces in complex spaces forms,
of which Theorem 2.1 is a particular case.

\begin{theorem} Let $\tilde{M}$ be a simply connected
Riemann surface with a conformal metric $\mathfrak{g}$ and
a holomorphic cubic differential $\Phi_0$ such that
\begin{equation}
K = c - 2 ||\Phi_0||^2,
\end{equation}
where $||\Phi_0||$ is the norm with respect to $\mathfrak{g}$ and $K$
is the Gaussian curvature of the metric.
Then there exists an $S^1$ family of isometric minimal Lagrangian
immersions $\, u_{\tau} : \tilde{M}  \to X(c)$
with the associated holomorphic cubic differentials $\, \Phi_{\tau} =
e^{i\tau} \, \Phi_0, \, e^{i\tau} \in S^1$.
The immersion is unique up to ambient isometry of $X(c)$ for each
holomorphic cubic differential $\Phi_{\tau}$.
This family exhausts all isometric minimal
Lagrangian immersion of $\, \tilde{M}$ into $\, X(c)$.
\end{theorem}
The theorem follows from an application of Frobenius theorem [Gr].
For the last part,
suppose $\Phi$ and $\Phi'$ are holomorphic cubic differentials
induced from two  isometric minimal Lagrangian immersions of a given
metric $\mathfrak{g}$. From  (23), the $0$-divisors of $\Phi$ and
$\Phi'$
are determined from the Gaussian curvature. Thus  $\Phi / \Phi'$ is a
holomorphic function of unit length, hence a constant.

The metric of a minimal Lagrangian surface in $X(c)$ admits a simple
intrinsic characterization.
\begin{theorem}
Suppose $\, \mathfrak{g}$ is a metric of a minimal Lagrangian surface
in a complex space form $\, X(c)$.
Then the associated metric
\begin{equation}
\tilde{\mathfrak{g}} = (\frac{c-K}{2})^{\frac{1}{3}} \, \mathfrak{g}
\notag
\end{equation}
is flat.
Conversely, let $\, \mathfrak{g}$ be a metric on a simply connected
surface $\tilde{M}$ such that the Gaussian curvature
\begin{equation}
K < c, \notag
\end{equation}
and that the associated metric $\, \tilde{\mathfrak{g}}$ is flat.
Then there exists an $\, S^1$ family of isometric minimal Lagrangian
immersions $\, u_{\tau} : \tilde{M}  \to X(c), \, e^{i\tau} \in
S^1$.
\end{theorem}
\emph{Proof.} $\;$ We present a proof for the case  $c=1$, and use the
notations in Section 2.
For the first part of the theorem, take a local coordinate $z$ on
$\tilde{M}$, away from the zero locus of $\Phi$, so that
$\Phi = h (\omega^1 + i \, \omega^2)^3 = (dz)^3$. Then the metric
\begin{equation}
( \, \Phi \circ \bar{\Phi} \, )^{\frac{1}{3}} = dz \circ d \bar{z} =
(h \, \bar{h})^{\frac{1}{3}} \mathfrak{g} \notag
\end{equation}
is flat. The converse follows by reverse of this argument and
Theorem 8.1.  $\square$

Given a minimal Lagrangian surface $u : M \to X(c)$, we now define
its associate family as follows.
Let $\pi : \tilde{M} \to M$ be the universal covering of $M$ with the
pulled back metric and holomorphic differential $\pi^*(\Phi)$.
By Theorem 8.1, there exist isometric minimal Lagrangian
immersions $u_{\tau} : \tilde{M} \to X(c)$
with the associated cubic differential $e^{i\tau} \pi^*(\Phi)$. A
surface $u_{\tau}(\tilde{M})$ is called an \emph{associate}
surface
of $u(M)$. $u_{\frac{\pi}{2}}(\tilde{M}) = M^{\star}$ is in
particular the \emph{conjugate} surface of $u(M)$.
Conjugate pair of surfaces satisfy the following dual reflection
principle.

Let  $\gamma$ be a curve on a minimal Lagrangian surface $M$ in
$X(c)$. $\gamma$  is \emph{Lagrangian} if
$\gamma \subset M \cap \mathcal{L}$ where $\mathcal{L}$ is a totally
geodesic Lagrangian surface in $X(c)$
that has orthogonal intersection with $M$. $\gamma \, $ is
\emph{complex} if $\gamma \subset M \cap \Pi$
where $\Pi$ is a totally geodesic complex curve in $X(c)$. In the
latter case $\Pi$ is necessarily orthogonal
to $M$. $\Pi$ ($\mathcal{L}$ respectively) is the fixed point locus of
an K\"ahler (anti- K\"ahler) involution of $X(c)$, and as remarked in Section 3,
the reflection principles continue to hold in this setting.

\begin{proposition}
A curve on a minimal Lagrangian surface in a complex space form   is
Lagrangian (complex respectively) if it is a geodesic
on which the associated holomorphic differential is real(purely
imaginary).
\end{proposition}
\emph{Proof.}$\;$ We present a proof for the case $c=1$. From the
structure equations (11), assume $e_2$ is tangent to the
curve $\gamma$, or $\omega^1=0$ on $\gamma$. Then the unique
anti-special Lagrangian 3-plane that contains $e_2 \wedge e_3 \,$ is
$\, - n_1 \wedge e_2 \wedge e_3$. A computation shows $d(n_1 \wedge
e_2 \wedge e_3)=0$ if and only if $\, \rho=0$ and $\, b=0$.
Similarly, the unique complex 2-plane
generated by $e_2 \wedge e_3$ is  $\, e_2 \wedge n_2 \wedge e _3
\wedge n_3$, and it is constant when $\rho=0$ and $\, a=0$. $\square$

\begin{corollary}
A Lagrangian geodesic on a minimal Lagrangian surface in a complex
space form corresponds to a complex geodesic
on its conjugate minimal Lagrangian surface, and vice versa.
\end{corollary}
For instance, if $M$ is obtained by successive reflection of a
fundamental domain across its Lagrangian or complex boundaries,
so is its conjugate $M^{\star}$.

\newpage

\begin{center} \noindent \textbf{BIBLIOGRAPHY} \end{center}

\noindent [Al] D. Allcock,  \emph{An Isoperimetric Inequality for
the Heisenberg Groups}, Geom. Funct. Anal. Vol 8, No 2 (1998), 219--233.

\noindent [AO] A. Amarzaya ; Y. Ohnita,  \emph{Hamiltonian stability
of certain minimal Lagrangian submanifolds in complex projective
spaces},  to appear in Tohoku Math. Journal.

\noindent [BCG] R. L. Bryant; S. S. Chern; R. B. Gardner; H. L.
Goldschmidt; P. A. Griffiths,  Exterior
differential systems. MSRI Publications, 18. Springer-Verlag, New
York, 1991.

\noindent [Br1] R. L. Bryant,  \emph{Conformal and minimal immersions
of compact surfaces into the $4$-sphere},
J. Diff. Geom. 17 (1982), no. 3, 455--473.

\noindent [Br2] $\underline{\qquad \quad}$, \emph{Submanifolds and
special structures on the octonians},  J. Differential Geom. 17
(1982), no. 2, 185--232.

\noindent [Cal] E. Calabi,  \emph{Complete affine hyperspheres I},
Symposia Mathematica, Vol. X,  pp. 19--38.

\noindent [Cox] H. S. M. Coxeter, \emph{Groups generated by unitary reflections of
period two}, Canad. Jour. Math. 9 (1957), 243--272

\noindent [Gr] P. A.  Griffiths, \emph{Some aspects of exterior
differential systems}, Complex geometry and Lie theory, 151--173,
Proc. Sympos. Pure Math., 53, Amer. Math. Soc., Providence, RI, 1991.

\noindent [Ha] M. Haskins, \emph{Special Lagrangian Cones},
arXiv: math. DG/0005164.

\noindent [HL1] R. Harvey ; H. B. Lawson Jr,  \emph{Calibrated
geometries},  Acta Math. 148 (1982), 47--157.

\noindent [HL2]  $\underline{\qquad \quad}$ ;  $\underline{\qquad
\quad}$, \emph{On boundaries of complex analytic varieties I}, Ann.
of Math. (2) 102 (1975), no. 2, 223--290,   $\,$ $II$. Ann. Math. (2)
106 (1977), no. 2, 213--238.

\noindent [Jo] D. Joyce, \emph{Special Lagrangian 3-folds and
integrable systems}, arXiv: math. DG /0101249.

\noindent [Ku]  Ernst Kuwert,  \emph{A compactness result for loops with
an $H^{\frac{1}{2}}$-bound}, J. reine angew. Math. 505 (1998) 1--22.

\noindent [La1]  H. B. Lawson Jr,  \emph{Complete minimal surfaces in
$S\sp{3}$}, Ann. of Math. (2) 92 (1970) 335--374.

\noindent [La2]  $\underline{\qquad \quad}$,  \emph{The global
behavior of minimal surfaces in $S\sp{n}$}, Ann. of Math. (2)
92 (1970) 224--237.

\noindent [LS]  $\underline{\qquad \quad}$ ; J. Simons,  \emph{On
stable currents and their application to global problems
in real and complex geometry}, Ann. of Math. (2) 98 (1973), 427--450.

\noindent [Ma]  Y. Matsuyama, \emph{Curvature pinching for totally
real submanifolds of a complex projective space},
J. Math. Soc. Japan 52 (2000), no. 1, 51--64.

\noindent [Mil] J. Milnor; J. Stasheff, Characteristic classes.
Annals of Mathematics Studies, No.76.
Princeton University Press, Princeton, N. J.; University of Tokyo
Press, Tokyo, 1974.

\noindent [MRU]  S. Montiel ; A. Ros ; F. Urbano,
\emph{Curvature pinching and eigenvalue rigidity for minimal
submanifolds}, Math. Z. 191 (1986), no. 4, 537--548.

\noindent [Na]  H. Naitoh, \emph{Totally real parallel submanifolds
in $P\sp{n}(c)$},
Tokyo J. Math. 4 (1981), no. 2, 279--306.

\noindent [Oh]  Y. Oh,  \emph{Second variation and stabilities of
minimal Lagrangian submanifolds in Kahler
manifolds}, Invent. Math. 101 (1990), no. 2, 501--519.

\noindent [Sal] S. Salamon, \emph{Harmonic and holomorphic maps},
161--224, Lecture Notes in Math., 1164, Springer, Berlin, 1985.

\noindent [ScU]  R. Schoen; K. Uhlenbeck,
\emph{Boundary regularity and the Dirichlet problem for harmonic maps},
 J. Diff. Geom. 18 (1983), 253--268.

\noindent [ScW]  $\underline{\qquad \quad}$; J. G. Wolfson,
\emph{Minimizing area among Lagrangian surfaces: the mapping problem},
arxiv: Math. DG/0008244.

\noindent [Si]  J. Simons, \emph{Minimal varieties in riemannian
manifolds}, Ann. of Math. (2) 88 (1968) 62--105.

\noindent [TeT] K. Tenenblat; C-L. Terng,  \emph{B\"acklund's theorem for n-dimensional
submanifolds of $\mathbb{R}^{2n-1}$},  Ann. of Math.111 (1980) 477--490.

\noindent [Wa1]  S. H. Wang, \emph{Calibrated lifts of minimal submanifolds},
arXiv. math. DG/0109214.

\noindent [Wa2]  S. H. Wang, \emph{Relative isoperimetric inequalities for Lagrangian
disks}, arXiv. math. DG/0301139.

\vspace{2pc}

\begin{flushleft}
\noindent Sung Ho Wang \\
Department of Mathematics \\
Postech\\
Pohang, Korea 790 - 784 \\
\emph{email}: \textbf{wang@postech.ac.kr}
\end{flushleft}

\end{document}